\documentclass[12pt]{amsart}
\usepackage{amsfonts}

\newtheorem{theorem}{Theorem}

\newtheorem{definition}[theorem]{Definition}
\newtheorem{example}[theorem]{Example}

\newtheorem{lemma}[theorem]{Lemma}

\newtheorem{proposition}[theorem]{Proposition}

\begin{document}
\title[Bernstein operators for exponential polynomials]{Bernstein operators for exponential polynomials}
\author[J.M. Aldaz, O. Kounchev and H. Render]{J.M. Aldaz, O. Kounchev and H. Render}
\thanks{The first and the last author are partially supported by Grant
MTM2006-13000-C03-03 of the D.G.I. of Spain. The last two authors
acknowledge support within the project ``Institutes Partnership'' with the
Alexander von Humboldt Foundation, Bonn.}
\address{J.M. Aldaz: Departamento de Matem\'{a}ticas y Computaci\'{o}n, Universidad
de La Rioja, Edificio Vives, Luis de Ulloa s/n., 26004 Logro\~{n}o, Spain.}
\email{aldaz@dmcunirioja.es}
\address{O. Kounchev: Institute of Mathematics and Informatics, Bulgarian Academy of
Sciences, 8 Acad. G. Bonchev Str., 1113 Sofia, Bulgaria.}
\email{kounchev@gmx.de}
\address{H. Render: Departamento de Matem\'{a}ticas y Computaci\'{o}n, Universidad de
La Rioja, Edificio Vives, Luis de Ulloa s/n., 26004 Logro\~{n}o, Espa\~{n}a.}
\email{render@gmx.de}

\begin{abstract}
Let $L$ be a linear differential operator with constant coefficients of
order $n$ and complex eigenvalues $\lambda _{0},...,\lambda _{n}$. Assume
that the set $U_{n}$ of all solutions of the equation $Lf=0$ is closed under
complex conjugation. If the length of the interval $\left[ a,b\right] $ is
smaller than $\pi /M_{n}$, where $M_{n}:=\max \left\{ \left| \text{Im}%
\lambda _{j}\right| :j=0,...,n\right\} $, then there exists a basis $p_{n,k}$%
, $k=0,...n$, of the space $U_{n}$ with the property that each $p_{n,k}$ has
a zero of order $k$ at $a$ and a zero of order $n-k$ at $b,$ and each $%
p_{n,k}$ is positive on the open interval $\left( a,b\right) .$ Under the
additional assumption that $\lambda _{0}$ and $\lambda _{1}$ are real and
distinct, our first main result states that there exist points $%
a=t_{0}<t_{1}<...<t_{n}=b$ and positive numbers $\alpha _{0},..,\alpha _{n}$%
, such that the operator
\begin{equation*}
B_{n}f:=\sum_{k=0}^{n}\alpha _{k}f\left( t_{k}\right) p_{n,k}\left( x\right)
\end{equation*}
satisfies $B_{n}e^{\lambda _{j}x}=e^{\lambda _{j}x}$, for $j=0,1.$ The
second main result gives a sufficient condition guaranteeing the uniform
convergence of $B_{n}f$ to $f$ for each $f\in C\left[ a,b\right] $.

2000 Mathematics Subject Classification: \emph{Primary: 41A35, Secondary
41A50}

Key words and phrases: \emph{Bernstein polynomial, Bernstein operator,
extended Chebyshev system, exponential polynomial}
\end{abstract}

\maketitle

\section{Introduction}

Let $\lambda _{0},...,\lambda _{n}$ be complex numbers, let $\Lambda _{n}$
be the vector $\left( \lambda _{0},...,\lambda _{n}\right) $, and de\-fine
the linear differential operator $L$ with constant coefficients by
\begin{equation}
L=\left( \frac{d}{dx}-\lambda _{0}\right) ....\left( \frac{d}{dx}-\lambda
_{n}\right) .  \label{eqdefL}
\end{equation}
Complex-valued solutions $f$ of the equation $Lf=0$ are called \emph{\
exponential polynomials} or $L$\emph{-polynomials.} They provide natural
generalizations of classical, trigonometric, and hyperbolic polynomials (see
\cite{Schu83}), and the so-called $\mathcal{D}$-polynomials considered in
\cite{MoNe00}. For example, it is well known that one can develop a nice
spline theory based on cardinal exponential polynomials (see e.g. \cite
{Micc76}, \cite{Schu81}, \cite{MaMe02}) and a satisfactory nonstationary
multiresolutional analysis for cardinal exponential splines, see the results
in \cite{dDR93}, \cite{oggybook}, \cite{kounchevrenderappr}, \cite
{kounchevrenderpams} and \cite{LySc93}, rediscovered in \cite{Unser}.
Another motivation stems from the investigation of a new class of
multivariate splines, the so-called polysplines, cf. \cite{oggybook}, \cite
{KounchevRenderJAT}.

Special interest in exponential polynomials has arisen recently within
Computer Aided Geometric Design for modelling parametric curves. On the one
hand, special systems of exponential polynomials are considered, such as
\begin{equation*}
1,x,...,x^{n-1},\cos x,\sin x,
\end{equation*}
(which corresponds to the case $\lambda _{0}=...=\lambda _{n-2}=0$, $\lambda
_{n-1}=i$, and $\lambda _{n}=-i),$ cf. \cite{Cost00}, \cite{MPS97}, \cite
{Zhan96} and \cite{CLM}, \cite{CMP07} for further generalizations. On the
other hand, a remarkable result is the existence of a so-called \emph{%
normalized} \emph{Bernstein basis} in certain classes of extended Chebyshev
systems, see \cite{CMP04}, \cite{Mazu05}. In order to explain this result,
let us recall that a subspace $U_{n}$ of $C^{n}\left( I\right) $, the space
of $n$-times continuously differentiable complex-valued functions on an
interval $I$, is called an \emph{extended Chebyshev system for the subset }$%
A\subset I$ if $U_{n}$ has dimension $n+1$ and each non-zero $f\in U_{n}$
vanishes at most $n$ times in $A$ (with multiplicities). A system $%
p_{n,k}\in U_{n},k=0,...,n$, is a \emph{Bernstein-like basis} for $a\neq
b\in I$, if the function $p_{n,k}$ has a zero of order $k$ at $a$, and a
zero of order $n-k$ at $b$ for $k=0,...,n$. For example, in the polynomial
case a Bernstein-like basis $P_{n,k}$ for $\left\{ a,b\right\} $ may be
defined explicitly by
\begin{equation}
P_{n,k}\left( x\right) :=\frac{1}{k!}\frac{1}{\left( b-a\right) ^{n-k}}%
\left( x-a\right) ^{k}\left( b-x\right) ^{n-k}.  \label{defpNk}
\end{equation}
The above-mentioned result in \cite{CMP04}, \cite{Mazu05} says the
following: Assume that the constant function $1$ is in $U_{n};$ clearly then
there exist coefficients $\alpha _{k},k=0,...,n,$ such that $%
1=\sum_{k=0}^{n}\alpha _{k}p_{n,k},$ since $p_{n,k},k=0,...,n$, is a basis.
The \emph{normalization property} proved in \cite{CMP04} and \cite{Mazu05}
for a certain class of Chebyshev systems says that the coefficients $\alpha
_{k}$ are \emph{positive}.

In this paper we shall be concerned with Bernstein-like bases and Bernstein
operators for the set of exponential polynomials induced by a linear
differential operator $L$ of the type (\ref{eqdefL}), i.e.
\begin{equation}
U_{n}=E_{\left( \lambda _{0},...,\lambda _{n}\right) }:=\left\{ f\in
C^{\infty }\left( \mathbb{R}\right) :Lf=0\right\} .  \label{Espace}
\end{equation}
It is easy to see that there exists a Bernstein-like basis $%
p_{n,k},k=0,...,n $ for $a\neq b$ if and only if $E_{\left( \lambda
_{0},...,\lambda _{n}\right) }$ is an extended Chebyshev system for the
\emph{set} $\left\{ a,b\right\} .$ In order to guarantee that the basis
functions $p_{n,k},k=0,...,n,$ are \emph{strictly positive} on the open
interval $\left( a,b\right) $ it is sufficient to know that $E_{\left(
\lambda _{0},...,\lambda _{n}\right) }$ is closed under complex conjugation
and that $E_{\left( \lambda _{0},...,\lambda _{n}\right) }$ is an extended
Chebyshev space for the closed \emph{interval} $\left[ a,b\right] .$ In
Section 2 we shall give the following criterion: $E_{\left( \lambda
_{0},...,\lambda _{n}\right) }$ is an extended Chebyshev space for the \emph{%
interval} $\left[ a,b\right] $ if $E_{\left( \lambda _{0},...,\lambda
_{n}\right) }$ is closed under complex conjugation and $b-a<\pi /M_{n}$,
where
\begin{equation}
M_{n}:=\max \left\{ \left| \text{Im}\lambda _{j}\right| :j=0,...,n\right\} .
\label{eqMM}
\end{equation}
Having at hand a Bernstein-like basis it is natural to ask whether one can
introduce a corresponding Bernstein operator, i.e. an operator of the type
\begin{equation}
B_{n}f\left( x\right) :=\sum_{k=0}^{n}\alpha _{k}f\left( t_{k}\right)
p_{n,k}\left( x\right)  \label{defBN}
\end{equation}
where the coefficients $\alpha _{0},...,\alpha _{n}$ and the knots $%
t_{0},...,t_{n}$ have to be defined in a suitable way. Our first main result
(Section 3, Theorem \ref{ThmBern}) states the following: assume that $%
E_{\left( \lambda _{0},...,\lambda _{n}\right) }$ is closed under complex
conjugation and $b-a<\pi /M_{n}.$ If $\lambda _{0}\neq \lambda _{1}$ are real%
\emph{\ }then there exist unique points $a=t_{0}<t_{1}<...<t_{n}=b$, and
unique \emph{positive} coefficients $\alpha _{0},...,\alpha _{n}$, such that
the operator $B_{n}:C\left[ a,b\right] \rightarrow E_{\left( \lambda
_{0},...,\lambda _{n}\right) }$ defined by (\ref{defBN}) has the following
reproducing property
\begin{equation}
B_{n}\left( e^{\lambda _{0}x}\right) =e^{\lambda _{0}x}\text{ and }%
B_{n}\left( e^{\lambda _{1}x}\right) =e^{\lambda _{1}x}.  \label{eqnorming}
\end{equation}
Of course, the latter property is reminiscent of the well known fact that
the classical (polynomial) Bernstein operator $B_{n}$ on $[0,1]$ satisfies $%
B_{n}1=1$ and $B_{n}x=x.$ Note that the assumption $b-a<\pi /M_{n}$ is
crucial: we give an example of an extended Chebyshev system $E_{\left(
\lambda _{0},...,\lambda _{n}\right) }$ over an \emph{interval} $\left[ a,b%
\right] $ (for $n=3)$ such that the condition (\ref{eqnorming}) implies the
positivity of the coefficients $\alpha _{0},\alpha _{1},\alpha _{2},\alpha
_{3}$ but the points $t_{0},t_{1},t_{2},t_{3}$ in (\ref{defBN}) are not
ordered, namely they satisfy the inequality $a=t_{0}<t_{2}<t_{1}<t_{3}=b.$
Additionally, we discuss the case $\lambda _{0}=\lambda _{1}$, the function $%
e^{\lambda _{1}x}$ in (\ref{eqnorming}) is replaced by $xe^{\lambda _{0}x}.$

It follows from the above construction that the operator $B_{n}$ defined by (%
\ref{defBN}) satisfying (\ref{eqnorming}) is a \emph{positive} operator.
Using a Korovkin-type theorem for extended Chebyshev systems we derive in
Section 4 a sufficient criterion for the uniform convergence of $B_{n}f$ on $%
f\in C\left[ a,b\right] .$ The criterion is formulated in terms of the basis
functions $p_{n,k}$ and their derivatives at the point $b $.

In Section 5 we consider exponential polynomials for equidistant eigenvalues
$\Lambda _{n}=\left( \lambda _{0},....,\lambda _{n}\right) $, i.e. for $%
\lambda _{j}=\lambda _{0}+j\omega $ for $j=0,...,n$. We briefly discuss the
relationship of a Bernstein-type theorem due to S. Morigi and M. Neamtu with
our results.

\section{Bernstein bases for complex eigenvalues}

In order to give the reader more intuition about exponential polynomials we
shall recall some elementary facts. In the case of pairwise different $%
\lambda _{j},j=0,...,n,$ the space $E_{\left( \lambda _{0},...,\lambda
_{n}\right) }$ is the linear span generated by the functions
\begin{equation*}
e^{\lambda _{0}x},e^{\lambda _{1}x},....,e^{\lambda _{n}x}.
\end{equation*}
When some $\lambda _{j}$ occurs $m_{j}$ times in $\Lambda _{n}=\left(
\lambda _{0},...,\lambda _{n}\right)$, a basis for the space $E_{\left(
\lambda _{0},...,\lambda _{n}\right) }$ is given by the linearly independent
functions
\begin{equation*}
x^{s}e^{\lambda _{j}x}\qquad \text{for }s=0,1,...,m_{j}-1.
\end{equation*}
We say that the vector $\Lambda _{n}\in \mathbb{C}^{n+1}$is \emph{equivalent}
to the vector $\Lambda _{n}^{\prime }\in \mathbb{C}^{n+1}$ if the
corresponding differential operators are equal (so the spaces of solutions
are equal). This is the same to say that each $\lambda $ occurs in $\Lambda
_{n}$ and $\Lambda _{n}^{\prime }$ with the same multiplicity. Since the
differential operator $L$ defined in (\ref{eqdefL}) does not depend on the
order of differentiation, it is clear that each permutation of the vector $%
\Lambda _{n}$ is equivalent to $\Lambda _{n}.$ Hence the space $E_{\left(
\lambda _{0},...,\lambda _{n}\right) }$ does not depend on the order of the
eigenvalues $\lambda _{0},...,\lambda _{n}.$

The $k$-th derivative of a function $f$ is denoted by $f^{\left( k\right) }.$
A function $f\in C^{n}\left( I,\mathbb{C}\right) $ has a \emph{zero of order
}$k$ or\emph{\ of multiplicity} $k$ at a point $a\in I$ if $f\left( a\right)
=...=f^{\left( k-1\right) }\left( a\right) =0$ and $f^{\left( k\right)
}\left( a\right) \neq 0.$ We shall repeatedly use the fact that
\begin{equation}
k!\cdot \lim_{x\rightarrow 0}\frac{f\left( x\right) }{\left( x-a\right) ^{k}}%
=f^{\left( k\right) }\left( a\right) .  \label{eqLim}
\end{equation}
for any function $f\in C^{\left( k\right) }(I)$ with $f\left( a\right)
=...=f^{\left( k-1\right) }\left( a\right) =0.$

\begin{definition}
A system of functions $p_{n,k},$ $k=0,...,n$ in the space $E_{\left( \lambda
_{0},...,\lambda _{n}\right) }$ is called \emph{Bernstein-like basis} of $%
E_{\left( \lambda _{0},...,\lambda _{n}\right) }$ for $a\neq b\in \mathbb{R}$
if and only if each function $p_{n,k}$ has a zero of order $k$ at $a$ and a
zero of order $n-k$ at $b$ for $k=0,...,n.$
\end{definition}

It is easy to see that a Bernstein-like basis $p_{n,k},$ $k=0,...,n$ (if it
exists) is indeed a basis for the space $E_{\left( \lambda _{0},...,\lambda
_{n}\right) }.$ Moreover the basis functions are unique up to a non-zero
multiplicative constant. In case of existence we shall require that
\begin{equation}
k!\lim_{x\rightarrow a,x>a}\frac{p_{n,k}\left( x\right) }{\left( x-a\right)
^{k}}=p_{n,k}^{\left( k\right) }\left( a\right) =1  \label{eqnorm}
\end{equation}
and we shall call $p_{n,k},$ $k=0,...,n$, \emph{the Bernstein basis} of $%
E_{\left( \lambda _{0},...,\lambda _{n}\right) }$ with respect to $a\neq b.$

In order to give a characterization of the existence of Bernstein bases, let
us recall the general fact (cf. \cite{Micc76}) that for $\Lambda _{n}=\left(
\lambda _{0},...,\lambda _{n}\right) \in \mathbb{C}^{n+1}$ there exists a
unique function $\Phi _{n}\in E_{\left( \lambda _{0},...,\lambda _{n}\right)
}$ such that $\Phi _{n}\left( 0\right) =....=\Phi _{n}^{\left( n-1\right)
}\left( 0\right) =0$ and $\Phi _{n}^{\left( n\right) }\left( 0\right) =1.$
An explicit formula for $\Phi _{n}$ is given by
\begin{equation}
\Phi _{n}\left( x\right) :=\Phi _{\Lambda _{n}}\left( x\right) :=\left[
\lambda _{0},...,\lambda _{n}\right] e^{xz}=\frac{1}{2\pi i}\int_{\Gamma
_{r}}\frac{e^{xz}}{\left( z-\lambda _{0}\right) ...\left( z-\lambda
_{n}\right) }dz  \label{defPhi}
\end{equation}
where $\left[ \lambda _{0},...,\lambda _{n}\right] $ denotes the divided
difference, and $\Gamma _{r}$ is a path in the complex plane defined by $%
\Gamma _{r}\left( t\right) =re^{it}$, $t\in \left[ 0,2\pi \right] $,
surrounding all the scalars $\lambda _{0},...,\lambda _{n}$. We shall call $%
\Phi _{n}$ the \emph{fundamental function}. Moreover we define
\begin{equation}
\Phi _{n,k}\left( x\right) :=\det \left(
\begin{array}{ccc}
\Phi _{n}\left( x\right) & ... & \Phi _{n}^{\left( k\right) }\left( x\right)
\\
\vdots &  & \vdots \\
\Phi _{n}^{\left( k\right) }\left( x\right) & ... & \Phi _{n}^{\left(
2k\right) }\left( x\right)
\end{array}
\right)  \label{defAnk}
\end{equation}
for each $k=0,...,n.$ The following characterization is straightforward, see
e.g. \cite{Veli07}:

\begin{theorem}
\label{ThmA}Let $\left( \lambda _{0},...,\lambda _{n}\right) \in \mathbb{C}%
^{n+1}$ and let $a\neq b\in \mathbb{R}.$ Then the following statements are
equivalent:

a) There exists a Bernstein basis $p_{n,k},$ $k=0,...,n$ in the space $%
E_{\left( \lambda _{0},...,\lambda _{n}\right) }$ for $\left\{ a,b\right\} ,$

b) $E_{\left( \lambda _{0},...,\lambda _{n}\right) }$ is an extended
Chebyshev system for $\left\{ a,b\right\} ,$

c) $\Phi _{n,k}\left( b-a\right) \neq 0$ for $k=0,...,n.$
\end{theorem}

The equivalence of a) and b) also holds in the context of Chebyshev systems
of real-valued functions, see \cite{CMP04}, \cite{GoMa03}, \cite{Mazu99},
\cite{Mazu05}. Further references on properties of Bernstein bases are \cite
{FaGo96} and \cite{Pena02}.

It is well known that for \emph{real} \emph{eigenvalues} $\lambda
_{0},...,\lambda _{n}$ the space $E_{\left( \lambda _{0},...,\lambda
_{n}\right) }$ is an extended Chebyshev system over \emph{any} interval $%
\left[ a,b\right] ,$ so a Bernstein basis exists in that case. In what
follows we want to discuss the case of complex eigenvalues. The reader
interested only in real eigenvalues may skip the rest of this section.

The following example is instructive:

\begin{example}
\label{Ex1}Let $\Lambda _{2}=\left( 0,i,-i\right) .$ Then $\Phi _{2}\left(
x\right) =1-\cos x$ is the fundamental function. Since $\Phi _{2}\left( 2\pi
k\right) =0$ it follows that $E_{\Lambda _{2}}$ does not possess a Bernstein
basis for $\left\{ 0,2\pi k\right\} $ with $k\in \mathbb{Z}.$ On the other
hand, if $b\neq 2\pi k,k\in \mathbb{Z},$ then $E_{\Lambda _{2}}$ does
possess a Bernstein basis for $\left\{ 0,b\right\} ,$ explicitly given by
\begin{eqnarray*}
p_{2,2}\left( x\right) &=&1-\cos x\text{, } \\
p_{2,1}\left( x\right) &=&\sin x-\sin b\frac{1-\cos x}{1-\cos b}, \\
p_{2,0}\left( x\right) &=&\frac{1-\cos \left( x-b\right) }{1-\cos b}.
\end{eqnarray*}
\end{example}

Note that $E_{\left( 0,i,-i\right) }$ possesses a Bernstein basis for $%
\left\{ 0,\pi \right\} $ but that the subspace $E_{\left( i,-i\right) }$
does not possess a Bernstein basis for $\left\{ 0,\pi \right\} $ since $%
\varphi _{\left( i,-i\right) }\left( x\right) =\sin x$ has then two zeros in
$\left\{ 0,\pi \right\} $.

Let us take now $b=3\pi $ in Example \ref{Ex1}. Then a Bernstein basis
exists for $\left\{ 0,3\pi \right\} $ but the basis function $p_{2,1}\left(
x\right) =\sin x$ takes negative values on the interval $\left[ 0,3\pi %
\right] .$ So Bernstein basis functions may fail to be positive.

Let us recall that $E_{\left( \lambda _{0},...,\lambda _{n}\right) }$ is
\emph{closed under complex conjugation} if for each $f\in $ $E_{\left(
\lambda _{0},...,\lambda _{n}\right) }$ the complex conjugate function $%
\overline{f}$ is again in $E_{\left( \lambda _{0},...,\lambda _{n}\right) }.$
It is easy to see that for complex numbers $\lambda _{0},...,\lambda _{n}$
the space $E_{\left( \lambda _{0},...,\lambda _{n}\right) }$ is closed under
complex conjugation if and only if there exists a permutation $\sigma $ of
the indices $\left\{ 0,...,n\right\} $ such that $\overline{\lambda _{j}}%
=\lambda _{\sigma \left( j\right) }$ for $j=0,...,n.$ In other words, $%
E_{\left( \lambda _{0},...,\lambda _{n}\right) }$ is closed under complex
conjugation if and only if the vector $\Lambda _{n}=\left( \lambda
_{0},...,\lambda _{n}\right) $ is equivalent to the conjugate vector $%
\overline{\Lambda _{n}}.$

\begin{proposition}
\label{PropCC}Suppose that $E_{\left( \lambda _{0},...,\lambda _{n}\right) }$
is an extended Chebyshev system for $a\neq b\in \mathbb{R}.$ Then the space $%
E_{\left( \lambda _{0},...,\lambda _{n}\right) }$ is closed under complex
conjugation if and only if the basis functions $p_{n,k}$ are real-valued on $%
\mathbb{R}$ for each $k=0,...,n.$
\end{proposition}

\begin{proof}
Suppose that $E_{\left( \lambda _{0},...,\lambda _{n}\right) }$ is closed
under complex conjugation. Then $\overline{p_{n,k}}$ is in $E_{\left(
\overline{\lambda _{0}},...,\overline{\lambda _{n}}\right) }=E_{\left(
\lambda _{0},...,\lambda _{n}\right) },$ it has a zero of order $k$ at $a $,
and a zero of order $n-k$ at $b.$ By uniqueness, $\overline{p_{n,k}}%
=Dp_{n,k} $ for some constant $D,$ and it follows from (\ref{eqnorm}) that $%
D=1.$ Thus $p_{\left( \lambda _{0},\lambda _{1},...,\lambda _{n}\right) ,k}$
is real-valued. The converse is easy, since the $p_{n,k}$ are real-valued
functions on $\mathbb{R,}$ and they form a basis.
\end{proof}

Assume $E_{\left( \lambda _{0},...,\lambda _{n}\right) }$ is closed under
complex conjugation and $E_{\left( \lambda _{0},...,\lambda _{n}\right) }$
is an extended Chebyshev system for the closed interval $\left[ a,b\right] $%
: then it easy to see that the Bernstein basis functions are strictly
positive on the interval $\left( a,b\right) $, and there exists $x_{0}\in
\left( a,b\right) $ such that $p_{n,k}$ is strictly increasing on $\left(
a,x_{0}\right) $ and decreasing on $\left( x_{0},b\right) ,$ see e.g. \cite
{Veli07}. Next we study when $E_{\left( \lambda _{0},...,\lambda _{n}\right)
}$ is an extended Chebyshev over an interval $\left[ a,b\right] .$

\begin{theorem}
Let $\left( \lambda _{0},...,\lambda _{n}\right) \in \mathbb{C}^{n+1}$ and
assume that $E_{\left( \lambda _{0},...,\lambda _{n}\right) }$ is closed
under complex conjugation. Then the following statements are equivalent for $%
a<b\in \mathbb{R}:$

a) $E_{\left( \lambda _{0},...,\lambda _{n}\right) }$ is an extended
Chebyshev system over the interval $\left[ a,b\right] ,$

b) $E_{\left( \lambda _{0},...,\lambda _{n}\right) }$ is an extended
Chebyshev system for all $\left\{ a,x\right\} $ with $x\in \left( a,b\right]
,$

c) The functions $x\longmapsto \Phi _{n,k}\left( x-a\right) ,$ $k=0,...,n,$
have no zeros in $\left( a,b\right] .$
\end{theorem}

\begin{proof}
Clearly $a)\rightarrow b)$ is trivial, and $b)$ and $c)$ are equivalent by
Theorem \ref{ThmA}. For $c)\rightarrow a)$ note that the function $\Phi _{n}$
is real-valued, so $\Phi _{n},\Phi _{n}^{\prime },...,\Phi _{n}^{\left(
n\right) }$ are real-valued and they form a basis of $E_{\left( \lambda
_{0},...,\lambda _{n}\right) }.$ We show that $f\in E_{\left( \lambda
_{0},...,\lambda _{n}\right) }$ has at most $n$ zeros in $\left[ a,b\right]
. $ Since $E_{\left( \lambda _{0},...,\lambda _{n}\right) }$ is closed under
complex conjugation we may assume that $f$ is real-valued. We can write $%
f=a_{0}\Phi _{n}+...+a_{n}\Phi _{n}^{\left( n\right) }.$ Then $f$ has a zero
of order $r\in \left\{ 0,..,n-1\right\} $ at $a,$ implying that $%
a_{n}=...=a_{n-r+1}=0.$ Hence $f$ is in the real linear span of $\Phi
_{n},...,\Phi _{n}^{\left( n-r\right) }$, which will be denoted by $U_{n-r}.$
Since $\Phi _{n,k}\left( b\right) \neq 0$ for $k=0,...,n-r$, there exists by
continuity some $\delta >0$ such that $\Phi _{n,k}\left( y\right) \neq 0$
for all $y\in \left[ b-\delta ,b+\delta \right] $ and $k=0,...,n-r.$ By
Theorem 2.3 in \cite[p. 52]{Karl68}, applied to the open interval $\left(
a,b+\delta \right) $, each function in $U_{r}$ has at most $n-r$ zeros
(counting the multiplicities) in the open interval $\left( a,b+\delta
\right) .$ Hence $f$ has at most $n$ zeros on $\left[ a,b\right] .$
\end{proof}

\begin{lemma}
\label{Lem0}If $E_{\left( \lambda _{0},....,\lambda _{n}\right) }$ is an
extended Chebyshev system over $\left[ a,b\right] $ and $\gamma $ is a real
number, then $E_{\left( \lambda _{0},....,\lambda _{n}\right) }$ is an
extended Chebyshev system over $\left[ a+\gamma ,b+\gamma \right] .$
\end{lemma}

\begin{lemma}
\label{Lem1}If $E_{\left( \lambda _{0},....,\lambda _{n}\right) }$ is an
extended Chebyshev system over $\left[ a,b\right] $ and $c$ is a complex
number then $E_{\left( \lambda _{0}-c,....,\lambda _{n}-c\right) }$ is an
extended Chebyshev system over $\left[ a,b\right] .$
\end{lemma}

\begin{proof}
If $f\in E\left( \lambda _{0},....,\lambda _{n}\right) $, then $g$ defined
by $g\left( x\right) =e^{-cx}f\left( x\right) $ is in $E_{\left( \lambda
_{0}-c,....,\lambda _{n}-c\right) }.$ If $g$ had more than $n$ zeros in $%
\left[ a,b\right] $ then $f$ would have more than $n$ zeros in $\left[ a,b%
\right] ,$ a contradiction.
\end{proof}

\begin{lemma}
\label{Lem2}If $E_{\left( \lambda _{0},....,\lambda _{n}\right) }$ is an
extended Chebyshev system over $\left[ a,b\right] $ and $c$ is a positive
number, then $E_{\left( \lambda _{0}c,....,\lambda _{n}c\right) }$ is an
extended Chebyshev system over $\left[ a,a+\frac{b-a}{c}\right] $.
\end{lemma}

\begin{proof}
By Lemma \ref{Lem0} we may assume that $a=0.$ If $f\in E_{\left( \lambda
_{0},....,\lambda _{n}\right) }$ then $g$, defined by $g\left( x\right)
:=f\left( cx\right) $, is in $E_{\left( \lambda _{0}c,....,\lambda
_{n}c\right) }.$ Suppose that $g$ has more than $n$ zeros in $\left[ 0,\frac{%
b}{c}\right] .$ Then $f$ has more than $n$ zeros in $\left[ 0,b\right] ,$ a
contradiction.
\end{proof}

The following is the main result of this section:

\begin{theorem}
\label{ThmB}Let $\left( \lambda _{0},...,\lambda _{n}\right) \in \mathbb{C}%
^{n+1}$ and assume that $E_{\left( \lambda _{0},...,\lambda _{n}\right) }$
is closed under complex conjugation. If $\left| \text{Im}\lambda _{j}\right|
\leq M_{n}$ for $j=0,...,n$, then $E_{\left( \lambda _{0},...,\lambda
_{n}\right) }$ is an extended Chebyshev system for the interval $\left[ a,b%
\right] $, provided $b-a<\pi /M_{n}.$
\end{theorem}

\begin{proof}
By an inductive argument, it suffices to prove the following two statements
for an extended Chebyshev system $E_{\left( \lambda _{0},....,\lambda
_{n}\right) }$ over $\left[ a,b\right] $, closed under complex conjugation:

1) If $\lambda _{n+1}$ is real then $E_{\left( \lambda _{0},....,\lambda
_{n+1}\right) }$ is an extended Chebyshev system for $\left[ a,b\right] ,$

2) If $\lambda _{n+1}$ is a non-real complex number, then $E\left( \lambda
_{0},....,\lambda _{n},\lambda _{n+1},\overline{\lambda _{n+1}}\right) $ is
an extended Chebyshev system over $\left[ a,d\right] $, for any $d$ with $%
a<d\leq b$ and $d-a<\frac{\pi }{\left| \text{Im}\lambda _{n+1}\right| } $.

For a proof of 1) we use a standard argument: let $f\in E_{\left( \lambda
_{0},....,\lambda _{n+1}\right) }$ be non-zero with $m$ zeros in $\left[ a,b%
\right] $. We may assume that $f$ is real-valued since $E_{\left( \lambda
_{0},....,\lambda _{n+1}\right) }$ is closed under complex conjugation. Then
$h\left( x\right) :=e^{-\lambda _{n+1}x}f\left( x\right) $ is real-valued
and it has $m$ zeros in $\left[ a,b\right] .$ By Rolle's theorem $h^{\prime
}\left( x\right) $ has at least $m-1$ zeros in $\left[ a,b\right] .$ Since
\begin{equation*}
e^{\lambda _{n+1}x}h^{\prime }\left( x\right) =e^{\lambda _{n+1}x}\frac{d}{dx%
}\left( e^{-\lambda _{n+1}x}f\left( x\right) \right) =\left( \frac{d}{dx}%
-\lambda _{n+1}\right) f\left( x\right) =:F\left( x\right)
\end{equation*}
we conclude that $F$ has at least $m-1$ zeros in $\left[ a,b\right] .$ But $%
F $ is in $E_{\left( \lambda _{0},....,\lambda _{n}\right) },$ so it has at
most $n$ zeros, and hence $m-1\leq n.$

For a proof of 2) note that by Lemma \ref{Lem1} we may assume that $c:=\text{%
Re}\lambda _{n+1}$ is zero. Without loss of generality let $\text{Im}\lambda
_{n+1}>0,$ and by Lemma \ref{Lem2} it suffices to prove 2) for the case that
$\lambda _{n+1}=i.$ It is clear that $E_{\left( \lambda _{0},....,\lambda
_{n},i,-i\right) }$ is closed under complex conjugation. Furthermore, by
Lemma \ref{Lem0} we may assume that $\left[ a,d\right] \subset I:=\left( -%
\frac{1}{2}\pi ,\frac{1}{2}\pi \right) .$ Let us introduce the auxiliary
function
\begin{equation*}
v\left( x\right) =\frac{1+\tan x}{1-\tan x}
\end{equation*}
defined on the interval $I.$ A computation shows that $v^{\prime }=v^{2}+1.$
Let $u$ be a primitive function of $v,$ so we have $u^{\prime }=v$ and $%
u^{\prime \prime }=\left( u^{\prime }\right) ^{2}+1.$ Let us define $%
g:=e^{u}.$ Then $g$ satisfies the differential equation $g^{\prime \prime
}g-2\left( g^{\prime }\right) ^{2}=g^{2}$ and a computation shows that
\begin{equation*}
g\frac{d}{dx}\left[ g^{-2}\frac{d}{dx}\left( gf\right) \right] =\left( \frac{%
d^{2}}{dx^{2}}+1\right) f.
\end{equation*}
Now we can argue as above: if $f\in E_{\left( \lambda _{0},....,\lambda
_{n},i,-i\right) }$ has $m$ zeros in the interval $\left[ a,d\right] $, so
does $gf.$ Thus $\frac{d}{dx}\left( gf\right) $ has at least $m-1$ zeros in $%
\left[ a,d\right] .$ Hence $g^{-2}\frac{d}{dx}\left( gf\right) $ has at
least $m-1$ zeros in $\left[ a,d\right] $ and we conclude that $\frac{d}{dx}%
\left[ g^{-2}\frac{d}{dx}\left( gf\right) \right] $ has at least $m-2$ zeros
in $\left[ a,d\right] $. Therefore $\left( \frac{d^{2}}{dx^{2}}+1\right) f$
has at least $m-2$ zeros in $\left[ a,d\right] $. Since $d\leq b$ and $%
\left( \frac{d^{2}}{dx^{2}}+1\right) f\in E_{\left( \lambda
_{0},....,\lambda _{n}\right) },$ and since $E_{\left( \lambda
_{0},....,\lambda _{n}\right) }$ is a Chebyshev system over $\left[ a,b%
\right] $ we obtain $m-2\leq n.$
\end{proof}

For a discussion of complex zeros of exponential polynomials we refer to the
recent work \cite{Wiel01}.

\section{Recursive relations for Bernstein bases}

Let $\left( \lambda _{0},...,\lambda _{n}\right) \in \mathbb{C}^{n+1}$ and
assume that $E_{\left( \lambda _{0},...,\lambda _{n}\right) }$ is an
extended Chebyshev system for $\left\{ a,b\right\} .$ We can construct a
Bernstein basis for $a\neq b\in \mathbb{R}$ via the following procedure: put
$q_{0}\left( x\right) =\Phi _{n}\left( x-a\right) ,$ which clearly has a
zero of order $n$ at $a.$ Then $q_{0}\left( b\right) \neq 0$ since $%
E_{\left( \lambda _{0},...,\lambda _{n}\right) }$ is an extended Chebyshev
system for $\left\{ a,b\right\} $. We define $q_{1}:=q_{0}^{\left( 1\right)
}-\alpha _{0}q_{0}$, where $\alpha _{0}=q_{0}^{\left( 1\right) }\left(
b\right) /q_{0}\left( b\right)$. Then $q_1$ has a zero of order $n-1$ at $a$
and a zero of order $1$ at $b.$ For $k\geq 2$ we define $q_{k}$ recursively
by
\begin{equation}
q_{k}:=q_{k-1}^{\left( 1\right) }-\left( \alpha _{k-1}-\alpha _{k-2}\right)
\cdot q_{k-1}-\beta _{k}q_{k-2}  \label{eqq}
\end{equation}
with coefficients $\alpha _{k-1},\alpha _{k-2}$ and $\beta _{k}$ to be
determined. Note that $q_{k}$ has a zero of order at least $k-2$ at $b,$ and
a zero of order $n-k$ at $a.$ The coefficients $\alpha _{k-1},\alpha _{k-2}$
and $\beta _{k}$ are chosen so that $q_{k}$ has a zero of order $k$ at $b,$
which is achieved by defining
\begin{equation}
\beta _{k}:=\frac{q_{k-1}^{\left( k-1\right) }\left( b\right) }{%
q_{k-2}^{\left( k-2\right) }\left( b\right) }\text{ and }\alpha _{k-1}:=%
\frac{q_{k-1}^{\left( k\right) }\left( b\right) }{q_{k-1}^{\left( k-1\right)
}\left( b\right) }.  \label{eqBP3}
\end{equation}
Then $p_{n,n-k}:=q_{k}$ for $k=0,...,n$ is the Bernstein basis satisfying
condition (\ref{eqnorm}).

The proof of the following proposition is easy and therefore omitted.

\begin{proposition}
\label{PropNull}Let $c\in \mathbb{C}$ and define $c+\Lambda _{n}:=\left(
c+\lambda _{0},....,c+\lambda _{n}\right) .$ If there exists a Bernstein
basis $p_{n,k},k=0,...,n$ for $E_{\Lambda _{n}}$ and $a\neq b$, then there
exists a Bernstein basis of $E_{c+\Lambda _{n}}$ given by
\begin{equation*}
p_{c+\Lambda _{n},k}\left( x\right) =p_{\Lambda _{n},k}\left( x\right)
e^{c\left( x-a\right) }
\end{equation*}
for $k=0,...,n.$
\end{proposition}

In the polynomial case the Bernstein basis $P_{n,k}$ defined in (\ref{defpNk}%
) satisfies the useful identity
\begin{equation*}
\frac{d}{dx}P_{n,k}=P_{n-1,k-1}-\frac{n-k}{b-a}P_{n-1,k},
\end{equation*}
which follows directly by differentiating (\ref{defpNk}). Next we present
its analog for exponential polynomials. In what follows we shall use the
more precise but lengthier notation $p_{\left( \lambda _{0},...,\lambda
_{n}\right) ,k}$ instead of $p_{n,k}$.

\begin{proposition}
\label{PropABL}Suppose that $E_{\left( \lambda _{0},...,\lambda _{n}\right)
} $ and $E_{\left( \lambda _{0},...,\lambda _{n-1}\right) }$ are extended
Chebyshev systems for $a\neq b\in \mathbb{R}.$ Define for $k=0,....,n-1$ the
numbers
\begin{equation}
d_{k}:=\lim_{x\uparrow b}\frac{\frac{d}{dx}p_{\left( \lambda
_{0},...,\lambda _{n}\right) ,k}\left( x\right) }{p_{\left( \lambda
_{0},...,\lambda _{n-1}\right) ,k}\left( x\right) }\neq 0.  \label{eqPR}
\end{equation}
Then,
\begin{equation}
\left( \frac{d}{dx}-\lambda _{n}\right) p_{\left( \lambda _{0},...,\lambda
_{n}\right) ,k}=p_{\left( \lambda _{0},...,\lambda _{n-1}\right)
,k-1}+d_{k}p_{\left( \lambda _{0},...,\lambda _{n-1}\right) ,k}
\label{eqPREC}
\end{equation}
for any $k=1,...,n-1.$ Furthermore, for $k=0$ the right hand side of (\ref
{eqPREC}) is equal to $d_{0}p_{\left( \lambda _{0},...,\lambda _{n-1}\right)
,0},$ while for $k=n$, it is equal to $p_{\left( \lambda _{0},...,\lambda
_{n-1}\right) ,n-1}.$
\end{proposition}

\begin{proof}
Let $f_{k}$ be the left hand side of (\ref{eqPREC}) and let $1\leq k\leq
n-1. $ Using the fact that $f_{k}$ has a zero of order $k-1$ at $a$ and a
zero of order $n-k-1$ at $b$, it is easy to see that $f_{k}=c_{k}p_{\left(
\lambda _{0},...,\lambda _{n-1}\right) ,k-1}+d_{k}p_{\left( \lambda
_{0},...,\lambda _{n-1}\right) ,k}$ for some constants $c_{k}$ and $d_{k}.$
Now simple limit considerations complete the proof.
\end{proof}

\begin{proposition}
\label{Propid} Suppose that $E_{\left( \lambda _{0},...,\lambda _{n}\right)
} $ and $E_{\left( \lambda _{0},...,\lambda _{n-1}\right) }$ are extended
Chebyshev systems for $a\neq b\in \mathbb{R.}$ Assume that $p_{\left(
\lambda _{0},\dots ,\lambda _{n-1},\lambda _{n}\right) ,k}\left( x\right)
=p_{\left( \lambda _{0},\dots ,\lambda _{n-1},\eta _{n}\right) ,k}\left(
x\right) $ on $\left( a,b\right) $ for a given $k$. Then $\lambda _{n}=\eta
_{n}.$ The same holds if instead of $\lambda_n$ we consider any other
eigenvalue.
\end{proposition}

\begin{proof}
Suppose first that $1\leq k\leq n-1$. By Proposition \ref{PropABL} we obtain
the equations
\begin{align*}
\left( \frac{d}{dx}-\lambda _{n}\right) p_{\left( \lambda _{0},\dots
,\lambda _{n-1},\lambda _{n}\right) ,k}& =p_{\left( \lambda _{0},...,\lambda
_{n-1}\right) ,k-1}+d_{k}p_{\left( \lambda _{0},...,\lambda _{n-1}\right)
,k}, \\
\left( \frac{d}{dx}-\eta _{n}\right) p_{\left( \lambda _{0},\dots ,\lambda
_{n-1},\eta _{n}\right) ,k}& =p_{\left( \lambda _{0},...,\lambda
_{n-1}\right) ,k-1}+D_{k}p_{\left( \lambda _{0},...,\lambda _{n-1}\right)
,k}.
\end{align*}
Using the order of the zeros at $b$ and assumption $p_{\left( \lambda
_{0},\dots ,\lambda _{n-1},\lambda _{n}\right) ,k}\left( x\right) =p_{\left(
\lambda _{0},\dots ,\lambda _{n-1},\eta _{n}\right) ,k}\left( x\right) $ we
get $d_{k}=D_{k},$ so the right hand sides are equal. Subtracting we
conclude that $\lambda _{n}=\eta _{n}.$ The cases $k=0$ and $k=1$
immediately follow from Proposition \ref{PropABL}.
\end{proof}

For the polynomial Bernstein basis over the interval $\left[ 0,1\right] $
one often uses identities like
\begin{equation}
1=\left( x+\left( 1-x\right) \right) ^{n}=\sum_{k=0}^{n}\frac{n!}{\left(
n-k\right) !}P_{n,k}\left( x\right)  \label{eq16}
\end{equation}
where $P_{n,k}$ is defined as in (\ref{defpNk}). The following is an analog
for exponential polynomials:

\begin{theorem}
\label{ThmRep}Suppose that $E_{\left( \lambda _{0},...,\lambda _{n}\right) }$
and $E_{\left( \lambda _{0},...,\lambda _{n-1}\right) }$ are extended
Chebyshev systems for $a\neq b\in \mathbb{R}.$ Let $d_{0},\dots ,d_{n-1}$ be
the non-zero numbers defined in (\ref{eqPR}). Then
\begin{equation}
e^{\left( x-a\right) \lambda _{n}}=p_{\left( \lambda _{0},\dots ,\lambda
_{n}\right) ,0}\left( x\right) +\sum_{k=1}^{n}\left( -1\right)
^{k}d_{0}\cdots d_{k-1}\cdot p_{\left( \lambda _{0},\dots ,\lambda
_{n}\right) ,k}\left( x\right) .  \label{pain2}
\end{equation}
Furthermore for $k=1,\dots ,n-1,$ we have the equality
\begin{equation}
d_{0}\cdots d_{k-1}=\left( -1\right) ^{n}\frac{e^{\left( b-a\right) \lambda
_{n}}}{p_{\left( \lambda _{0},...,\lambda _{n}\right) ,n}\left( b\right) }%
\frac{1}{d_{k}\cdots d_{n-1}}.  \label{eqlambda2}
\end{equation}
\end{theorem}

\begin{proof}
Write $e^{\left( x-a\right) \lambda _{n}}=\sum_{k=0}^{n}\beta _{k}p_{\left(
\lambda _{0},\dots ,\lambda _{n}\right) ,k}\left( x\right) $ for
coefficients $\beta _{0},...,\beta _{n}$. Inserting $x=a$ yields $1=\beta
_{0}p_{\left( \lambda _{0},...,\lambda _{n}\right) ,0}\left( 0\right) ,$ so $%
\beta _{0}=1$ by (\ref{eqnorm}). Proposition \ref{PropABL} yields
\begin{align*}
& 0=\left( \frac{d}{dx}-\lambda _{n}\right) e^{\left( x-a\right) \lambda
_{n}}=d_{0}p_{\left( \lambda _{0},...,\lambda _{n-1}\right) ,0}(x)+\beta
_{1}p_{\left( \lambda _{0},...,\lambda _{n-1}\right) ,0}(x) \\
& +\sum_{k=1}^{n-2}\left( \beta _{k+1}+\beta _{k}d_{k}\right) p_{\left(
\lambda _{0},...,\lambda _{n-1}\right) ,k}(x)+p_{\left( \lambda
_{0},...,\lambda _{n-1}\right) ,n-1}(x)\left[ \beta _{n-1}d_{n-1}+\beta _{n}%
\right] .
\end{align*}
Thus $\beta _{1}=-d_{0},$ $\beta _{k+1}=-\beta _{k}d_{k}$ for $k=1,...,n-2,$
and $\beta _{n-1}d_{n-1}+\beta _{n}=0.$ Hence, for $k=1,\dots ,n$ we have $%
\beta _{k}=(-1)^{k}d_{0}\cdots d_{k-1},$ and then (\ref{pain2}) follows.
Furthermore, by inserting $x=b$ in (\ref{pain2}) and recalling that $%
p_{\left( \lambda _{0},...,\lambda _{n}\right) ,k}$ has a zero of order $n-k$
at $x=b,$ we see that
\begin{equation*}
\left( -1\right) ^{n}d_{0}\cdots d_{n-1}\cdot p_{\left( \lambda _{0},\dots
,\lambda _{n}\right) ,n}\left( b\right) =e^{\left( b-a\right) \lambda _{n}}.
\end{equation*}
Thus, we get (\ref{eqlambda2}).
\end{proof}

Theorem \ref{ThmRep} does not hold when the assumption of having an extended
Chebyshev system $E_{\left( \lambda _{0},...,\lambda _{n-1}\right) }$ for $%
a\neq b$ is dropped: in Example \ref{Ex1}, with $b=\pi$, one has that
\begin{equation*}
1=\frac{1}{2}\left( 1-\cos x\right) +0\cdot \sin x+\frac{1}{2}\left( 1+\cos
x\right) ,
\end{equation*}
so $1$ is a linear combination of the Bernstein basis functions with a \emph{%
zero} coefficient.

\begin{theorem}
\label{ThmMono}Suppose that $E_{\left( \lambda _{0},...,\lambda _{n}\right)
} $, $E_{\left( \lambda _{0},\lambda _{2}...,\lambda _{n}\right) }$, $%
E_{\left( \lambda _{1}...,\lambda _{n}\right) }$ and $E_{\left( \lambda
_{2}...,\lambda _{n}\right) }$ are extended Chebyshev systems for $a\neq
b\in \mathbb{R}.$ Let $\lambda _{0}\neq \lambda _{1}.$ Then there exists a
constant $C_{k}^{\lambda _{0},\lambda _{1}}\left( \Lambda _{n}\right) \neq 0$
such that
\begin{equation}
p_{\left( \lambda _{0},\lambda _{2},...,\lambda _{n}\right) ,k}-p_{\left(
\lambda _{1},\lambda _{2},...,\lambda _{n}\right) ,k}=C_{k}^{\lambda
_{0},\lambda _{1}}\left( \Lambda _{n}\right) \cdot p_{\left( \lambda
_{0},\lambda _{1},\lambda _{2},...,\lambda _{n}\right) ,k+1}.  \label{eqB}
\end{equation}
Moreover,
\begin{equation}
\lim_{x\rightarrow b}\frac{p_{\left( \lambda _{0},\lambda _{2},...,\lambda
_{n}\right) ,k}\left( x\right) }{p_{\left( \lambda _{1},...,\lambda
_{n}\right) ,k}\left( x\right) }\neq 1.  \label{eqB2}
\end{equation}
\end{theorem}

\begin{proof}
Let $B\left( x\right) $ be the function on the left hand side of (\ref{eqB}%
). Then $B$ has a zero of order $k+1$ at $a,$ since $p_{\left( \lambda
_{0},\lambda _{2},...,\lambda _{n}\right) ,k}$ and $p_{\left( \lambda
_{1},\lambda _{2},...,\lambda _{n}\right) ,k}$ have a zero of order $k$ at $%
a $, and $B^{\left( k\right) }\left( a\right) =\lim_{x\rightarrow a}\frac{%
B\left( x\right) }{\left( x-a\right) ^{k}}=0$ by (\ref{eqnorm}). Furthermore
$B$ has a zero of order $n-k-1$ at $b.$ By Proposition \ref{Propid}, $B$ is
not identically zero (here we need that $E_{\left( \lambda _{2}...,\lambda
_{n}\right) }$ is an extended Chebyshev system for $a\neq b)$. Since $B\in
E_{\left( \lambda _{0},\lambda _{1},...,\lambda _{n}\right) },$ it must be a
non-zero multiple of $p_{\left( \lambda _{0},\lambda _{1},...,\lambda
_{n}\right) ,k+1}$.

Finally, suppose that the limit in (\ref{eqB2}) is equal to $1.$ Then $%
p_{\left( \lambda _{0},\lambda _{2},...,\lambda _{n}\right) ,k}^{\left(
n-1-k\right) }\left( b\right) =p_{\left( \lambda _{1},...,\lambda
_{n}\right) ,k}^{\left( n-1-k\right) }\left( b\right) .$ By (\ref{eqB}) we
conclude that $p_{\left( \lambda _{0},\lambda _{1}\lambda _{2},...,\lambda
_{n}\right) ,k+1}^{\left( n-1-k\right) }\left( b\right) =0.$ Hence $%
p_{\left( \lambda _{0},\lambda _{1}\lambda _{2},...,\lambda _{n}\right)
,k+1} $ has a zero of order $n-k$ at $b$ and a zero of order $k+1$ at $a,$ a
contradiction.
\end{proof}

In the case of equidistant eigenvalues it is possible to define a Bernstein
basis explicitly:

\begin{proposition}
Suppose that $\omega \neq 0$ and $\lambda _{j}=\lambda _{0}+j\omega $ for $%
j=0,...,n.$ Then
\begin{equation}
p_{n,k}\left( x\right) :=\frac{e^{\lambda _{0}\left( x-a\right) }}{k!\omega
^{k}}\left( e^{\omega \left( x-a\right) }-1\right) ^{k}\left( \frac{%
1-e^{\omega \left( x-b\right) }}{1-e^{\omega \left( a-b\right) }}\right)
^{n-k}  \label{eqdefppp}
\end{equation}
is a Bernstein basis for $\left[ a,b\right] $ satisfying
\begin{equation}
k!\lim_{x\rightarrow a}p_{n,k}\left( a\right) /\left( x-a\right) ^{k}=1.
\label{normalizati3}
\end{equation}
\end{proposition}

\begin{proof}
It is easy to see that $p_{n,k}$ is an exponential polynomial. Furthermore, $%
p_{n,k}$ has a zero at $x=a$ of order $k$ and a zero of order $n-k$ at $x=b.$
\end{proof}

\begin{lemma}
\label{LemC} Let $\omega \neq 0$ and let $\lambda _{j}=\lambda _{0}+j\omega $
for $j=0,...,n$. Then
\begin{equation*}
p_{\left( \lambda _{0},\lambda _{1},...,\lambda _{n-1}\right) ,k}-p_{\left(
\lambda _{1},\lambda _{2},...,\lambda _{n}\right) ,k}=-\left( k+1\right)
\omega p_{\left( \lambda _{0},\lambda _{1}\lambda _{2},...,\lambda
_{n}\right) ,k+1}.
\end{equation*}
\end{lemma}

\begin{proof}
This is a computation using (\ref{eqdefppp}).
\end{proof}

The following result is crucial for the proof of the existence of a
Bernstein operator. In this theorem we shall use a homotopy argument for the
eigenvalues $\Lambda _{n}$ and the assumption (\ref{assumM}) will guarantee
that the corresponding Bernstein bases with respect to the points $a\neq b$
exist.

\begin{theorem}
\label{ThmMono2} Suppose that $E_{\left( \lambda _{2},...,\lambda
_{n}\right) }$ is closed under complex conjugation and $0<b-a<\pi /M_{n}$,
where
\begin{equation}
M_{n}=\max \left\{ \left| \text{Im}\lambda _{j}\right| :\text{for }%
j=2,...,n\right\} \text{.}  \label{assumM}
\end{equation}
If $\lambda _{0},\lambda _{1}\in \mathbb{R}$ and $\lambda _{0}<\lambda _{1}$%
, then
\begin{equation}
\lim_{x\rightarrow b}\frac{p_{\left( \lambda _{0},\lambda _{2},...,\lambda
_{n}\right) ,k}\left( x\right) }{p_{\left( \lambda _{1},...,\lambda
_{n}\right) ,k}\left( x\right) }<1.  \label{eqkl1}
\end{equation}
Furthermore, the function of $\lambda \in \mathbb{R}$
\begin{equation}
\lambda \longmapsto p_{\left( \lambda ,\lambda _{2},...,\lambda _{n}\right)
,k}\left( x\right)  \label{eqff}
\end{equation}
is strictly increasing for each $x\in \left( a,b\right) $.
\end{theorem}

\begin{proof}
By Proposition \ref{PropCC}, $p_{\left( \lambda ,\lambda _{2},...,\lambda
_{n}\right) ,k}\left( x\right) $ is real valued for every real $\lambda $.
Let now $\lambda _{0}<\lambda _{1}$ be real. It follows that $C_{k}^{\lambda
_{0},\lambda _{1}}\left( \Lambda _{n}\right) $ in (\ref{eqB}) is real.
Clearly, for $\lambda _{0}<\lambda _{1}$ and fixed $\left( \lambda
_{2},...,\lambda _{n}\right) $ the function in (\ref{eqff}) is increasing if
$C_{k}^{\lambda _{0},\lambda _{1}}\left( \Lambda _{n}\right) $ in (\ref{eqB}%
) is negative. From the inductive formula (\ref{eqq}) we get that the
function $\left( \lambda _{0},...,\lambda _{n}\right) \longmapsto p_{\left(
\lambda _{0},...,\lambda _{n}\right) ,k}$ is continuous. Now (\ref{eqB})
implies that $\left( \lambda _{0},...,\lambda _{n}\right) \longmapsto
C_{k}^{\lambda _{0},\lambda _{1}}\left( \Lambda _{n}\right) $ is continuous.
For $\lambda _{0}<\lambda _{1}$, define $\mu _{j}:=\lambda _{0}+j\left(
\lambda _{1}-\lambda _{0}\right) /n$, where $j=0,...,n.$ Since $\mu
_{0}=\lambda _{0}$ and $\mu _{n}=\lambda _{1}$ the function
\begin{equation*}
t\longmapsto C_{k}^{t\lambda _{0}+\left( 1-t\right) \mu _{0},t\lambda
_{1}+\left( 1-t\right) \mu _{n}}\left( t\left( \lambda _{0},...,\lambda
_{n}\right) +\left( 1-t\right) \left( \mu _{0}...,\mu _{n}\right) \right)
\end{equation*}
is continuous, and by Theorem \ref{ThmMono}, it has no zero on $\mathbb{R}.$
It follows that this function must have constant sign. Hence it suffices to
show that $C_{k}^{\mu _{0},\mu _{n}}(\mu _{0},...,\mu _{n})<0.$ But this
follows from Lemma \ref{LemC}. Thus, $\lambda \longmapsto p_{\left( \lambda
,\lambda _{2},...,\lambda _{n}\right) ,k}\left( x\right) $ is increasing, so
for $a<x<b$ and $\lambda _{0}<\lambda _{1}$,
\begin{equation}
p_{\left( \lambda _{0},\lambda _{2},...,\lambda _{n}\right) ,k}\left(
x\right) <p_{\left( \lambda _{1},\lambda _{2},...,\lambda _{n}\right)
,k}\left( x\right) .  \label{eigench}
\end{equation}
Dividing (\ref{eigench}) by its right hand side, using (\ref{eqB2}), and
taking the limit $x\uparrow b$, we get (\ref{eqkl1}).
\end{proof}

\begin{example}
\label{Ex2}Let us take $\lambda \in \mathbb{R}$. Then the fundamental
function $\Phi _{\left( \lambda ,i,-i\right) }$ with respect to $\left(
\lambda ,i,-i\right) $ is given by
\begin{equation*}
\Phi _{\left( \lambda ,i,-i\right) }\left( x\right) =\frac{e^{\lambda
x}-\cos x-\lambda \sin x}{\lambda ^{2}+1}.
\end{equation*}
Since $\Phi _{\left( \lambda ,i,-i\right) }$ is equal to the basis function $%
p_{\left( \lambda ,i,-i\right) ,2}$ it follows that $\lambda \longmapsto
\Phi _{\left( \lambda ,i,-i\right) }\left( x\right) $ is increasing for any $%
x$ in the interval $\left( 0,\pi \right) .$ Differentiating, it is easy to
check that $\lambda \longmapsto \Phi _{\lambda }\left( x\right) $ is
decreasing whenever $x<0$ is small in absolute value, and $\lambda $ is
sufficiently large. For $x$ with $\pi <x<2\pi $ it might be checked that $%
\lambda \longmapsto \Phi _{\lambda }\left( x\right) $ is not increasing.
\end{example}

\section{Construction of the Bernstein operator}

We now proceed to our first main result which roughly says the following:
given two functions $e^{\lambda _{0}x}$ and $e^{\lambda _{1}x}$ in the
extended Chebyshev system $E_{\left( \lambda _{0},...,\lambda _{n}\right) }$
over $\left[ a,b\right] $ we can find points $t_{0},...,t_{n}$ in the
interval $\left[ a,b\right] $ and positive numbers $\alpha _{0},...,\alpha
_{n}$ such that the operator $B_{n}$ defined by (\ref{defBB}) below
reproduces (or preserves) the functions $e^{\lambda _{0}x}$ and $e^{\lambda
_{1}x},$ i.e. (\ref{defsatB}) holds.

\begin{theorem}
\label{ThmBern} Let $\lambda _{0},...,\lambda _{n}$ be complex numbers with $%
\lambda _{0}$ and $\lambda _{1}$ real and $\lambda _{0}<\lambda _{1}.$
Suppose $E_{\left( \lambda _{0},...,\lambda _{n}\right) }$ is closed under
complex conjugation and $0<b-a<\pi /M_{n}$, where
\begin{equation*}
M_{n}:=\max \left\{ \left| \text{Im}\lambda _{j}\right| :j=0,...,n\right\} .
\end{equation*}
Define inductively points $t_{0},...,t_{n}$ by setting $t_{0}=a$ and
\begin{equation*}
e^{\left( \lambda _{0}-\lambda _{1}\right) \left( t_{k}-t_{k-1}\right)
}=\lim_{x\rightarrow b}\frac{p_{\left( \lambda _{0},\lambda _{2},...,\lambda
_{n}\right) ,k-1}\left( x\right) }{p_{\left( \lambda _{1},...,\lambda
_{n}\right) ,k-1}\left( x\right)}
\end{equation*}
for $k=1,2,...,n.$ Then
\begin{equation*}
a=t_{0}<t_{1}<....<t_{n}=b.
\end{equation*}
Put $\alpha _{0}=1,$ and define numbers
\begin{equation}
\alpha _{k}=e^{-\lambda _{0}\left( t_{k}-a\right) }\left( -1\right)
^{k}\prod_{l=0}^{k-1}\lim_{x\rightarrow b}\frac{\frac{d}{dx}p_{\left(
\lambda _{0},...,\lambda _{n}\right) ,l}\left( x\right) }{p_{\left( \lambda
_{1},...,\lambda _{n}\right) ,l}\left( x\right) }  \label{eqdefak}
\end{equation}
for $k=1,....,n.$ Then $\alpha _{0},...,\alpha _{n}>0$ and the operator $%
B_{\left( \lambda _{0},...,\lambda _{n}\right) }$ on $C\left[ a,b\right] $
defined by
\begin{equation}
B_{\left( \lambda _{0},...,\lambda _{n}\right) }f=\sum_{k=0}^{n}\alpha
_{k}f\left( t_{k}\right) p_{\left( \lambda _{0},...,\lambda _{n}\right) ,k}
\label{defBB}
\end{equation}
fixes the functions $e^{\lambda _{0}x}$ and $e^{\lambda _{1}x},$ i.e.
\begin{equation}
B_{\left( \lambda _{0},...,\lambda _{n}\right) }\left( e^{\lambda
_{0}x}\right) =e^{\lambda _{0}x}\text{ and }B_{\left( \lambda
_{0},...,\lambda _{n}\right) }\left( e^{\lambda _{1}x}\right) =e^{\lambda
_{1}x}.  \label{defsatB}
\end{equation}
Moreover, the real numbers $t_{0},...,t_{n}$ and $\alpha _{0},...,\alpha
_{n} $ satisfying (\ref{defsatB}) are unique.
\end{theorem}

\begin{proof}
First let us only assume that $E_{\left( \lambda _{0},\lambda
_{1},...,\lambda _{n}\right) },$ $E_{\left( \lambda _{1},\lambda
_{2},...,\lambda _{n}\right) }$ and $E_{\left( \lambda _{0},\lambda
_{2}...,\lambda _{n}\right) }$ are extended Chebyshev systems over $\left\{
a,b\right\} $, in order to clarify where we need the assumption $0<b-a<\pi
/M_{n}.$ Let $\beta _{0},....,\beta _{n}$ and $\gamma _{0},...,\gamma _{n}$
be the unique \emph{non-zero} coefficients, found in Theorem \ref{ThmRep},
that satisfy
\begin{equation}
e^{\lambda _{0}\left( x-a\right) }=\sum_{k=0}^{n}\beta _{k}p_{\left( \lambda
_{0},...,\lambda _{n}\right) ,k}\left( x\right) \text{ and }e^{\lambda
_{1}\left( x-a\right) }=\sum_{k=0}^{n}\gamma _{k}p_{\left( \lambda
_{0},...,\lambda _{n}\right) ,k}\left( x\right) .  \label{eqeq}
\end{equation}
The reproducing property of the Bernstein operator for $e^{\lambda _{0}x}$
in (\ref{defsatB}) implies that
\begin{equation*}
\sum_{k=0}^{n}e^{\lambda _{0}\left( t_{k}-a\right) }\alpha _{k}p_{\left(
\lambda _{0},...,\lambda _{n}\right) ,k}\left( x\right) =\sum_{k=0}^{n}\beta
_{k}p_{\left( \lambda _{0},...,\lambda _{n}\right) ,k}\left( x\right) .
\end{equation*}
Since $p_{\left( \lambda _{0},...,\lambda _{n}\right) ,k}$ is a basis we
conclude that $e^{\lambda _{0}\left( t_{k}-a\right) }\alpha _{k}=\beta _{k}.$
Similarly, $e^{\lambda _{1}\left( t_{k}-a\right) }\alpha _{k}=\gamma _{k}$
follows from $B_{\left( \lambda _{0},...,\lambda _{n}\right) }\left(
e^{\lambda _{1}x}\right) =e^{\lambda _{1}x}$. Dividing, we see that $t_{k}$
satisfies the equation
\begin{equation}
e^{\left( \lambda _{0}-\lambda _{1}\right) t_{k}}=\frac{\beta _{k}}{\gamma
_{k}}e^{\left( \lambda _{0}-\lambda _{1}\right) a}.  \label{eqtk}
\end{equation}
Hence, for $\alpha _{k}$ we obtain
\begin{equation}
\alpha _{k}=e^{-\lambda _{0}\left( t_{k}-a\right) }\beta _{k}.  \label{eqak}
\end{equation}

It is easy to see that $\beta _{k}$ and $\gamma _{k}$ are real, since the
functions $p_{\left( \lambda _{0},...,\lambda _{n}\right) ,k}$ are
real-valued (Proposition \ref{PropCC}) and both $\lambda _{0}$ and $\lambda
_{1}$ are real. Now $t\longmapsto e^{\left( \lambda _{0}-\lambda _{1}\right)
t}$ is real-valued and injective, so $t_{k}$ is \emph{uniquely} determined
by (\ref{eqtk}), and hence so is $\alpha _{k}$ by (\ref{eqak}). Moreover, it
is easy to see that the points $t_{0}=a$ and $t_{n}=b$ satisfy (\ref{eqtk}).

Next we want to show that $\alpha _{k}$ is positive: by Theorem \ref{ThmRep}
(applied to $\lambda _{0}$ instead of $\lambda _{n}$) we have $\beta
_{k}=\left( -1\right) ^{k}\widetilde{d}_{0}...\widetilde{d}_{k-1}$, where $%
\widetilde{d}_{l}$ is given by
\begin{equation}
\widetilde{d}_{l}=\lim_{x\rightarrow b}\frac{\frac{d}{dx}p_{\left( \lambda
_{0},...,\lambda _{n}\right) ,l}\left( x\right) }{p_{\left( \lambda
_{1},...,\lambda _{n}\right) ,l}\left( x\right) }.  \label{deka}
\end{equation}
The positivity of $p_{\left( \lambda _{1},...,\lambda _{n}\right) ,k}\left(
x\right) $ on $\left( a,b\right) $ implies that $\widetilde{d}_{l}$ is
negative. Hence, for $k=1,...,n-1$, equation (\ref{eqak}) yields
\begin{equation}
\alpha _{k}=e^{-\lambda _{0}\left( t_{k}-a\right) }\beta _{k}=e^{-\lambda
_{0}\left( t_{k}-a\right) }\left( -1\right) ^{k}\widetilde{d}_{0}...%
\widetilde{d}_{k-1}  \label{eqalphak}
\end{equation}
showing that $\alpha _{k}$ is positive. Similarly, $\gamma _{k}=\left(
-1\right) ^{k}D_{0}....D_{k-1}$ where $D_{l}$ is given by
\begin{equation}
D_{l}=\lim_{x\rightarrow b}\frac{\frac{d}{dx}p_{\left( \lambda
_{0},...,\lambda _{n}\right) ,l}\left( x\right) }{p_{\left( \lambda
_{0},\lambda _{2},...,\lambda _{n}\right) ,l}\left( x\right) }.  \label{Deka}
\end{equation}
Thus the points $t_{k}$ are defined by
\begin{equation}
e^{\left( \lambda _{0}-\lambda _{1}\right) t_{k}}=\frac{\beta _{k}}{\gamma
_{k}}e^{\left( \lambda _{0}-\lambda _{1}\right) a}=\frac{\widetilde{d}_{0}...%
\widetilde{d}_{k-1}}{D_{0}....D_{k-1}}e^{\left( \lambda _{0}-\lambda
_{1}\right) a}.  \label{eqdefft}
\end{equation}
Note that for $k=1,...,n$
\begin{equation}
e^{\left( \lambda _{0}-\lambda _{1}\right) \left( t_{k}-t_{k-1}\right) }=%
\frac{e^{\left( \lambda _{0}-\lambda _{1}\right) t_{k}}}{e^{\left( \lambda
_{0}-\lambda _{1}\right) t_{k-1}}}=\frac{\widetilde{d}_{k-1}}{D_{k-1}}%
=\lim_{x\rightarrow b}\frac{p_{\left( \lambda _{0},\lambda _{2},...,\lambda
_{n}\right) ,k-1}\left( x\right) }{p_{\left( \lambda _{1},...,\lambda
_{n}\right) ,k-1}\left( x\right) }.  \label{eqtkmono}
\end{equation}
Next we show that $t_{k}$ is in the interval $\left[ a,b\right] .$ Since $%
t_{0}=a$ and $t_{n}=b$, it suffices to show that $t_{k-1}<t_{k}.$ Since $%
\lambda _{0}<\lambda _{1}$ the requirement $t_{k-1}<t_{k}$ is equivalent to
the requirement that
\begin{equation*}
\lim_{x\rightarrow b}\frac{p_{\left( \lambda _{0},\lambda _{2},...,\lambda
_{n}\right) ,k-1}\left( x\right) }{p_{\left( \lambda _{1},...,\lambda
_{n}\right) ,k-1}\left( x\right) }<1.
\end{equation*}
Theorem \ref{ThmMono2} tells us this is true under the assumption that $%
\left| b-a\right| <\pi /M_{n}$, thus finishing the proof.
\end{proof}

In Example \ref{Ex2} we have computed the fundamental function for $\left(
\lambda ,i,-i\right) $ and it is easy to see that the Bernstein basis
function $p_{\left( \lambda ,i,-i\right) ,1}$ for $\left\{ 0,b\right\} $ is
given by
\begin{equation}
p_{\left( \lambda ,i,-i\right) ,1}\left( x\right) =\frac{\lambda e^{\lambda
x}+\sin x-\lambda \cos x}{\lambda ^{2}+1}-\frac{\left( e^{\lambda x}-\cos
x-\lambda \sin x\right) \left( \lambda e^{\lambda b}+\sin b-\lambda \cos
b\right) }{\left( e^{\lambda b}-\cos b-\lambda \sin b\right) \left( \lambda
^{2}+1\right) }.  \label{eqpp}
\end{equation}
Simple computations show that $p_{\left( \lambda ,i,-i\right) ,1}^{\prime
}\left( b\right) =\frac{e^{\lambda b}\cos b-\lambda e^{\lambda b}\sin b-1}{%
e^{\lambda b}-\cos b-\lambda \sin b}$ and
\begin{equation}
\frac{p_{\left( -\lambda ,i,-i\right) ,1}^{\prime }\left( b\right) }{%
p_{\left( \lambda ,i,-i\right) ,1}^{\prime }\left( b\right) }=\left( \frac{%
e^{\lambda b}-\cos b-\lambda \sin b}{1-e^{\lambda b}\cos b+\lambda
e^{\lambda b}\sin b}\right) ^{2}.  \label{eqquot}
\end{equation}
Consider now the Bernstein operator for $\Lambda _{3}=\left(
-1,1,i,-i\right) $ for the interval $\left[ 0,3.5\right] .$ Using Theorem
\ref{ThmA} it can be seen that $E_{\left( 1,i,-i\right) }$ and $E_{\left(
-1,i,-i\right) }$ are extended Chebyshev systems, at least for the interval $%
\left[ 0,3.8\right] .$ By property 1) in the proof of Theorem \ref{ThmB}, $%
E_{\left( -1,1,i,-i\right) }$ is an extended Chebyshev system for $\left[
0,3.8\right] .$ Hence by the proof of Theorem \ref{ThmRep} it follows from (%
\ref{eqtkmono}) that
\begin{equation*}
e^{-2\left( t_{2}-t_{1}\right) }=\frac{p_{\left( -1,i,-i\right) ,1}^{\prime
}\left( 3.5\right) }{p_{\left( 1,i,-i\right) ,1}^{\prime }\left( 3.5\right) }%
\approx 2.\,8454>1,
\end{equation*}
so $t_{2}-t_{1}$ must be \emph{negative}, and thus the Bernstein operator
for $\Lambda _{3}=\left( -1,1,i,-i\right) $ has the property that $%
t_{0}<t_{2}<t_{1}<t_{3}.$

If we take $b=\pi $ then (\ref{eqpp}) shows that $p_{\left( \lambda
,i,-i\right) ,1}=\sin x.$ Thus Theorem \ref{ThmMono} is not valid if we drop
the assumption that $E_{\left( \lambda _{2}...,\lambda _{n}\right) }$ is an
extended Chebyshev system for $a\neq b\in \mathbb{R}.$

By a limiting process one can handle the case that the eigenvalues $\lambda
_{0},\lambda _{1}$ are equal when replacing $e^{\lambda _{1}x}$ by $%
xe^{\lambda _{0}x}$in (\ref{defsatB}). However, we have not been able to
show that in this case the nodes $a=t_{0}\leq t_{1}\leq ....\leq t_{n}=b$
are distinct.

\begin{theorem}
Let $\lambda _{0},...,\lambda _{n}$ be complex numbers such that $\lambda
_{0}=\lambda _{1}$ is real. Suppose that $E_{\left( \lambda _{0},...,\lambda
_{n}\right) }$ is closed under complex conjugation and that $0<b-a<\pi
/M_{n} $ for
\begin{equation*}
M_{n}:=\max \left\{ \left| \text{Im}\lambda _{j}\right| :j=0,...,n\right\} .
\end{equation*}
Then there exist unique nodes $a=t_{0}\leq t_{1}\leq ....\leq t_{n}=b$ and
unique positive numbers $\alpha _{0},...,\alpha _{n}$ such that the operator
defined for $f\in C\left[ a,b\right] $ by
\begin{equation}
B_{\left( \lambda _{0},...,\lambda _{n}\right) }f=\sum_{k=0}^{n}\alpha
_{k}f\left( t_{k}\right) p_{\left( \lambda _{0},...,\lambda _{n}\right) ,k}
\label{defBern2}
\end{equation}
fixes the functions $e^{\lambda _{0}x}$ and $xe^{\lambda _{0}x}$.
\end{theorem}

\begin{proof}
For $\varepsilon \geq 0$ define $\Lambda _{\varepsilon }:=\left( \lambda
_{0},\lambda _{0}+\varepsilon ,\lambda _{2},...,\lambda _{n}\right) .$ By
Theorem \ref{ThmBern} there exist for each $\varepsilon >0$ points $%
a=t_{0}\left( \varepsilon \right) <t_{1}\left( \varepsilon \right)
<....<t_{n}\left( \varepsilon \right) =b$ and positive numbers $\alpha
_{0}\left( \varepsilon \right) ,...,\alpha _{n}\left( \varepsilon \right) $
such that the corresponding Bernstein operator $B_{\Lambda _{\varepsilon }}$
satisfies
\begin{equation*}
B_{\Lambda _{\varepsilon }}\left( e^{\lambda _{0}x}\right) =e^{\lambda _{0}x}%
\text{ and }B_{\Lambda _{\varepsilon }}\left( e^{\left( \lambda
_{0}+\varepsilon \right) x}\right) =e^{\left( \lambda _{0}+\varepsilon
\right) x}.
\end{equation*}
By compactness of the interval $\left[ a,b\right] $ there exists a sequence
of positive numbers $\varepsilon _{m}\rightarrow 0$ such that $t_{j}\left(
\varepsilon _{m}\right) $ converges to numbers $t_{j}$ for $m\rightarrow
\infty $ and for each $j=0,...,n.$ Clearly one has $a=t_{0}\leq t_{1}\leq
....\leq t_{n}=b.$ Let us write $e^{\lambda _{0}\left( x-a\right)
}=\sum_{k=0}^{n}\beta _{k}\left( \varepsilon \right) p_{\Lambda
_{\varepsilon },k}\left( x\right) $ for $\varepsilon \geq 0.$ Clearly $%
p_{\Lambda _{\varepsilon _{m}},k}\left( x\right) $ converges to $p_{\Lambda
_{0},k}\left( x\right) $ for $m\rightarrow \infty $ and $k=0,...,n,$ and $%
\beta _{k}\left( \varepsilon _{m}\right) $ converges to $\beta _{k}\left(
0\right) $ (cf. formula (\ref{eqPR}) and Theorem \ref{ThmRep}). By Theorem
\ref{ThmRep} $\beta _{k}\left( 0\right) $ is positive and the formula $%
e^{\lambda _{0}\left( t_{k}\left( \varepsilon \right) -a\right) }\alpha
_{k}\left( \varepsilon \right) =\beta _{k}\left( \varepsilon \right) $ now
shows that $\alpha _{k}\left( \varepsilon _{m}\right) $ converges to the
\emph{positive} number $\alpha _{k}\left( 0\right) .$ We define now the
Bernstein operator $B_{\Lambda _{0}}f$ by (\ref{defBern2}) with $\alpha
_{k}:=\alpha _{k}\left( 0\right) $ for $k=0,...,n.$ It is easy to see that
\begin{equation*}
B_{\Lambda _{0}}\left( e^{\lambda _{0}x}\right) =\lim_{\varepsilon
\rightarrow 0}B_{\Lambda _{\varepsilon }}\left( e^{\lambda _{0}x}\right)
=e^{\lambda _{0}x}.
\end{equation*}
Clearly
\begin{equation*}
f_{m}\left( x\right) :=\frac{e^{\left( \lambda _{0}+\varepsilon _{m}\right)
x}-e^{\lambda _{0}x}}{\varepsilon _{m}}\rightarrow xe^{\lambda _{0}x}\text{
for }m\rightarrow \infty
\end{equation*}
and since $B_{\Lambda _{\varepsilon _{m}}}f_{m}\left( x\right) =f_{m},$ a
limit argument shows that $B_{\Lambda _{0}}\left( xe^{\lambda _{0}x}\right)
=xe^{\lambda _{0}x}.$ The uniqueness is proven in a similar way as in the
last proof.
\end{proof}

\section{Convergence of the Bernstein operator}

Next we present a sufficient condition for the Bernstein operator $B_{\left(
\lambda_{0},...,\lambda_{n}\right) }$ to converge to the identity.

\begin{definition}
For each $n\in \mathbb{N}$, let $\{a\left( n,k\right) :k=0,....,n\}$ be a
triangular array of complex numbers. We say that $a\left( n,k\right) $
converges uniformly to $c$ if for each $\varepsilon >0$ there exists a
natural number $n_{0}$ such that $\left| a\left( n,k\right) -c\right|
<\varepsilon $, for all $n\geq n_{0}$ and all $k=0,...,n.$
\end{definition}

The following lemma is implicitly contained in \cite[p. 47]{Lore86}. For
completeness we include the proof.

\begin{lemma}
\label{PropLore}Let $\gamma >0.$ For each $n\in\mathbb{N}$ and each $%
j=0,....,n$, let $a\left( n,j\right)\in \left( 0,1\right)$ and $b\left(
n,j\right)\in\mathbb{R}$. Suppose that
\begin{equation}
\lim_{n\rightarrow \infty }\frac{\log b\left( n,j\right) }{\log a\left(
n,j\right) }=\gamma >0,  \label{eqdlog}
\end{equation}
and assume that the convergence is uniform in $j.$ Define $A_{k}\left(
n\right) =\prod_{j=k}^{n}a\left( n,j\right) $ and $B_{k}\left( n\right)
=\prod_{j=k}^{n}b\left( n,j\right) .$ Then $\lim_{n\rightarrow \infty
}\left( A_{k}\left( n\right) ^{\gamma }-B_{k}\left( n\right) \right) =0$
uniformly in $k.$
\end{lemma}

\begin{proof}
We have to show that for each $\varepsilon _{1}>0$ there exists an $n_{0}$
such that for all $n\geq n_{0}$ and all $k=0...,n,$
\begin{equation*}
\left| A_{k}\left( n\right) ^{\gamma }-B_{k}\left( n\right) \right|
<\varepsilon _{1}.
\end{equation*}
Fix $\varepsilon _{1}>0$, and select $\varepsilon \in (0,1)$ such that $%
1-\varepsilon ^{\varepsilon }<\varepsilon _{1}$, $\varepsilon +\varepsilon
^{\gamma }<\varepsilon _{1},$ $\varepsilon <\gamma $, and $\varepsilon
^{\gamma -\varepsilon }<\varepsilon _{1}.$ By (\ref{eqdlog}), there exists
an $n_{0}$ such that if $n\geq n_{0}$ and $j=0,...,n$,
\begin{equation*}
\left| \frac{\log b\left( n,j\right) }{\log a\left( n,j\right) }-\gamma
\right| <\varepsilon .
\end{equation*}
Then $\gamma -\varepsilon <\frac{\log b\left( n,j\right) }{\log a\left(
n,j\right) }<\gamma +\varepsilon $. Observe that $\log b\left( n,j\right) <0$
since $\log a\left( n,j\right) <0.$ So $\left( \gamma -\varepsilon \right)
\log a\left( n,j\right) >\log b\left( n,j\right) $ and $\left( \gamma
+\varepsilon \right) \log a\left( n,j\right) <\log b\left( n,j\right) .$
Hence $a\left( n,j\right) ^{\gamma +\varepsilon }\leq b\left( n,j\right) $
and $b\left( n,j\right) \leq a\left( n,j\right) ^{\gamma -\varepsilon }.$
Thus we have proven that
\begin{equation}
A_{k}\left( n\right) ^{\gamma +\varepsilon }\leq B_{k}\left( n\right) \text{
and }B_{k}\left( n\right) \leq A_{k}\left( n\right) ^{\gamma -\varepsilon }.
\label{eqaaakkn}
\end{equation}
Next we consider two cases: First assume that $A_{k}(n)^{\gamma }\geq
B_{k}\left( n\right) .$ Then, using (\ref{eqaaakkn}),
\begin{equation*}
0\leq A_{k}\left( n\right) ^{\gamma }-B_{k}\left( n\right) \leq A_{k}\left(
n\right) ^{\gamma }-A_{k}\left( n\right) ^{\gamma +\varepsilon
}=A_{k}(n)^{\gamma }\left( 1-A_{k}\left( n\right) ^{\varepsilon }\right) .
\end{equation*}
If $A_{k}\left( n\right) \geq \varepsilon $, then for all $n\geq n_{0}$ and
all $k$ we have, using that $A_{k}\left( n\right) <1$,
\begin{equation*}
0\leq A_{k}\left( n\right) ^{\gamma }-B_{k}\left( n\right) \leq
1-A_{k}\left( n\right) ^{\varepsilon }\leq 1-\varepsilon ^{\varepsilon
}<\varepsilon _{1}.
\end{equation*}
If $A_{k}\left( n\right) <\varepsilon $ for some $k,n$, then $B_{k}\left(
n\right) \leq A_{k}\left( n\right) ^{\gamma }\leq \varepsilon ^{\gamma }$
and $0\leq A_{k}\left( n\right) ^{\gamma }-B_{k}\left( n\right) \leq
\varepsilon +\varepsilon ^{\gamma }<\varepsilon _{1}.$In the second case we
have $A_{k}(n)^{\gamma }\leq B_{k}\left( n\right) .$ Then, from (\ref
{eqaaakkn}) we get
\begin{equation*}
0\leq B_{k}\left( n\right) -A_{k}\left( n\right) ^{\gamma }\leq A_{k}\left(
n\right) ^{\gamma -\varepsilon }-A_{k}\left( n\right) ^{\gamma
}=A_{k}(n)^{\gamma -\varepsilon }\left( 1-A_{k}\left( n\right) ^{\varepsilon
}\right) .
\end{equation*}
If $A_{k}\left( n\right) \geq \varepsilon $ we obtain
\begin{equation*}
0\leq B_{k}\left( n\right) -A_{k}\left( n\right) ^{\gamma }\leq
1-A_{k}\left( n\right) ^{\varepsilon }\leq 1-\varepsilon ^{\varepsilon
}<\varepsilon _{1}.
\end{equation*}
Finally, if $A_{k}\left( n\right) <\varepsilon $ for some $k,n$, then $0\leq
B_{k}\left( n\right) -A_{k}\left( n\right) ^{\gamma }\leq A_{k}\left(
n\right) ^{\gamma -\varepsilon }-A_{k}\left( n\right) ^{\gamma }\leq
\varepsilon ^{\gamma -\varepsilon }<\varepsilon _{1}.$
\end{proof}

Next we present our second main result:

\begin{theorem}
\label{ThmCon}Let $\lambda _{0},\lambda _{1},\lambda _{2}$ be distinct real
numbers and let $\Lambda _{n}=\left( \lambda _{0},\lambda _{1},....,\lambda
_{n}\right)$, where for $j=3,...,n$ the complex numbers $\lambda _{j}$ are
allowed to vary. Suppose each $E_{\left( \lambda _{0},...,\lambda
_{n}\right) }$ is closed under complex conjugation, and furthermore there
exists a positive number $M$ such that for every $n\ge 2$ and every $%
j=0,...,n$, we have $\left| \text{Im}\lambda _{j}\right| \leq M$. For each $%
k\leq n$ set
\begin{align}
a\left( n,k\right) & :=\lim_{x\rightarrow b}\frac{p_{\left( \lambda
_{0},\lambda _{2},....,\lambda _{n}\right) ,k}\left( x\right) }{p_{\left(
\lambda _{1},\lambda _{2},....,\lambda _{n}\right) ,k}\left( x\right) },%
\text{ and }  \label{assump1} \\
b\left( n,k\right) & :=\lim_{x\rightarrow b}\frac{p_{\left( \lambda
_{0},\lambda _{1},\lambda _{3},....,\lambda _{n}\right) ,k}\left( x\right) }{%
p_{\left( \lambda _{1},\lambda _{2},....,\lambda _{n}\right) ,k}\left(
x\right) }.  \label{assump2}
\end{align}
Let $t_{k}$, $k=0,\dots,n$, be the uniquely determined points given by
Theorem \ref{ThmBern}. Assume that
\begin{equation}
\lim_{n\to\infty} t_{k} -t_{k-1}= 0  \label{tktk1}
\end{equation}
uniformly in $k$, and likewise, that
\begin{equation}
\lim_{n\to\infty}\frac{\log b\left( n,k\right) }{t_{k}-t_{k+1}}=
\lambda_{2}-\lambda _{0}  \label{assumptonb}
\end{equation}
uniformly in $k$. Then the Bernstein operator $B_{\left( \lambda
_{0},...,\lambda _{n}\right) } $ defined in Theorem \ref{ThmBern}, converges
to the identity operator on $C([a,b],\mathbb{C})$ with the uniform norm.
\end{theorem}

\begin{proof}
1. We remind the reader that $\widetilde{d}_{k}$ is given by (\ref{deka}),
and $D_{k}$ by (\ref{Deka}). Recall that $B_{\left(
\lambda_{0},...,\lambda_{n}\right) }f\left( x\right) =\sum_{k=0}^{n}f\left(
t_{k}\right) \alpha_{k}p_{\left( \lambda_{0},...,\lambda_{n}\right)
,k}\left( x\right)$, where $\alpha _{k}=e^{-\lambda_{0}\left( t_{k}-a\right)
}\left( -1\right) ^{k}\widetilde{d}_{0}\cdots\widetilde{d}_{k-1}.$ By
construction we have $B_{\left( \lambda_{0},...,\lambda_{n}\right)
}e^{\lambda_{j}x}=e^{\lambda_{j}x}$ for $j=0,1.$ If we show that $B_{\left(
\lambda_{0},...,\lambda_{n}\right) }e^{\lambda_{2}\left( x-a\right) }$
converges to $e^{\lambda_{2}\left( x-a\right) },$ then it follows from the
generalized Korovkin theorem for Chebyshev systems that $B_{\Lambda_{n}}$
converges to the identity operator (cf. \cite{Koro60}, Theorem 8).

2. We may assume that $\lambda _{2}>0.$ Indeed, we can always translate the
vector $\Lambda _{n}=\left( \lambda _{0},...,\lambda _{n}\right) $ by a
positive constant $c$ so that $\lambda _{2}+c>0$. Then $E_{\left( \lambda
_{0}+c,...,\lambda _{n}+c\right) }$ is again closed under complex
conjugation; by Proposition \ref{PropNull}, the corresponding numbers in (%
\ref{assump1}) and (\ref{assump2}) are the same for the translated vector $%
c+\Lambda _{n}$. If the Bernstein operator converges for $c+\Lambda _{n}$,
then so does $B_{\left( \lambda _{0},...,\lambda _{n}\right) }$. Thus, we
may assume that $\lambda _{2}>0$ to begin with.

3. Write (by Theorem \ref{ThmRep}) $e^{\left( x-a\right) \lambda
_{2}}=\sum_{k=0}^{n}q_{k}p_{\left( \lambda _{0},...,\lambda _{n}\right)
,k}\left( x\right) $, where $q_{k}=\left( -1\right) ^{k}\widetilde{%
\widetilde{d}}_{0}\cdots \widetilde{\widetilde{d}}_{k-1}$ and
\begin{equation}
\widetilde{\widetilde{d}}_{n,k}:=\lim_{x\rightarrow b}\frac{\frac{d}{dx}%
p_{\left( \lambda _{0},...,\lambda _{n}\right) ,k}\left( x\right) }{%
p_{\left( \lambda _{0},\lambda _{1},\lambda _{3},....,\lambda _{n}\right)
,k}\left( x\right) }.  \label{detilde}
\end{equation}
Thus, for $\varphi \left( x\right) :=e^{\lambda _{2}\left( x-b\right)
}=e^{\lambda _{2}\left( x-a\right) }e^{\lambda _{2}\left( a-b\right) }$, by (%
\ref{eqalphak}) we have
\begin{equation}
B_{\left( \lambda _{0},...,\lambda _{n}\right) }\varphi (x)-\varphi
(x)=\sum_{k=0}^{n}\left( e^{\lambda _{2}\left( t_{k}-b\right) }-\frac{q_{k}}{%
\alpha _{k}}e^{\lambda _{2}\left( a-b\right) }\right) \alpha _{k}p_{\left(
\lambda _{0},...,\lambda _{n}\right) ,k}\left( x\right) .  \label{diffe}
\end{equation}
Also, $e^{\lambda _{2}\left( t_{k}-b\right) }=\prod_{j=k}^{n-1}e^{\lambda
_{2}\left( t_{j}-t_{j+1}\right) }$ for $k\leq n-1$. Define
\begin{equation}
\widetilde{a}\left( n,j\right) :=e^{t_{j}-t_{j+1}}<1.  \label{atilde}
\end{equation}
Since $\varphi \left( x\right) =e^{\lambda _{2}\left( x-a\right) }e^{\lambda
_{2}\left( a-b\right) }$, we have
\begin{equation}
\frac{q_{k}}{\alpha _{k}}e^{\lambda _{2}\left( a-b\right) }=\frac{\widetilde{%
\widetilde{d}}_{0}....\widetilde{\widetilde{d}}_{k-1}}{\widetilde{d}_{0}....%
\widetilde{d}_{k-1}}e^{\lambda _{0}\left( t_{k}-a\right) }e^{\lambda
_{2}\left( a-b\right) }.  \label{ddtilde}
\end{equation}
By (\ref{eqlambda2}) (applied to $\lambda _{0}$ instead of $\lambda _{n}$ )
we obtain
\begin{equation}
\frac{q_{k}}{\alpha _{k}}e^{\lambda _{2}\left( a-b\right) }=\frac{\widetilde{%
d}_{k}....\widetilde{d}_{n-1}}{\widetilde{\widetilde{d}}_{k}....\widetilde{%
\widetilde{d}}_{n-1}}e^{\lambda _{0}\left( t_{k}-b\right)
}=\prod_{j=k}^{n-1}\left( \frac{\widetilde{d_{j}}}{\widetilde{\widetilde{d}%
_{j}}}e^{\lambda _{0}\left( t_{j}-t_{j+1}\right) }\right) .  \label{identity}
\end{equation}
Since the left hand side of (\ref{identity}) is real for each $k=0,...,n-1$,
it is clear that
\begin{equation}
\widetilde{b}\left( n,j\right) :=\frac{\widetilde{d_{j}}}{\widetilde{%
\widetilde{d}_{j}}}e^{\lambda _{0}\left( t_{j}-t_{j+1}\right) }=b\left(
n,j\right) e^{\lambda _{0}\left( t_{j}-t_{j+1}\right) }  \label{btilde}
\end{equation}
is real for all $j=0,...,n-1$ and all $n.$ Suppose we know that
\begin{equation}
\frac{\log \widetilde{b}\left( n,j\right) }{\log \widetilde{a}\left(
n,j\right) }\rightarrow \lambda _{2}.  \label{eqbschlange}
\end{equation}
Then we may use Lemma \ref{PropLore}: Set $\widetilde{A}_{k}\left( n\right)
=\prod_{j=k}^{n-1}\widetilde{a}\left( n,j\right) $ and $\widetilde{B}%
_{k}\left( n\right) =\prod_{j=k}^{n-1}\widetilde{b}\left( n,j\right) .$ By
uniform convergence, given $\varepsilon >0$ there exists $n_{0}$ such that $%
\left| \widetilde{A}_{k}\left( n\right) ^{\lambda _{2}}-\widetilde{B}%
_{k}\left( n\right) \right| <\varepsilon $ for all $n\geq n_{0}$ and all $%
k\leq n$. Note that by (\ref{btilde}) and (\ref{identity}),
\begin{equation}
\widetilde{B}_{k}\left( n\right) =\frac{q_{k}}{\alpha _{k}}e^{\lambda
_{2}\left( a-b\right) }.  \label{painf}
\end{equation}
Since $\widetilde{A}_{k}\left( n\right) ^{\lambda _{2}}=e^{\lambda
_{2}\left( t_{k}-b\right) }$ and $1\leq e^{\lambda _{2}\left( t_{k}-a\right)
}$, from (\ref{diffe}) we get
\begin{equation}
\left| B_{\left( \lambda _{0},...,\lambda _{n}\right) }\varphi (x)-\varphi
(x)\right| \leq \varepsilon \sum_{k=0}^{n}\alpha _{k}p_{\left( \lambda
_{0},...,\lambda _{n}\right) ,k}\left( x\right)  \label{estimate}
\end{equation}
\begin{equation*}
\leq \varepsilon \sum_{k=0}^{n}e^{\lambda _{2}\left( t_{k}-a\right) }\alpha
_{k}p_{\left( \lambda _{0},...,\lambda _{n}\right) ,k}\left( x\right)
=\varepsilon e^{\lambda _{2}\left( b-a\right) }B_{\left( \lambda
_{0},...,\lambda _{n}\right) }\varphi (x)
\end{equation*}
for all $n\geq n_{0}$ and all $k\leq n$. So for every $x\in \lbrack a,b]$,
\begin{equation*}
\frac{1}{1+\varepsilon e^{\lambda _{2}\left( b-a\right) }}\varphi \left(
x\right) \leq B_{\left( \lambda _{0},...,\lambda _{n}\right) }\varphi \left(
x\right) \leq \frac{1}{1-\varepsilon e^{\lambda _{2}\left( b-a\right) }}%
\varphi \left( x\right) ,
\end{equation*}
proving uniform convergence of $B_{\left( \lambda _{0},...,\lambda
_{n}\right) }\varphi $ to to $\varphi$ on $\left[ a,b\right]$.

We show next that (\ref{eqbschlange}) holds. From formula (\ref{eqdefft}) we
obtain
\begin{equation}
a\left( n,k\right) =\frac{\widetilde{d}_{k}}{D_{k}}=e^{\left( \lambda
_{0}-\lambda_{1}\right) \left( t_{k+1}-t_{k}\right) }.  \label{asin}
\end{equation}

By assumption (\ref{tktk1}), $\lim_{n\longrightarrow\infty}a\left(
n,k\right) =1$ uniformly in $k.$ Since $\lambda_{0}\neq\lambda_{1}$, this
implies that

\begin{equation}
\lim_{n\to\infty}\frac{1-b\left( n,k\right) }{1-a\left( n,k\right) }=\frac{{%
\lambda _{2}-\lambda _{0}}}{{\lambda _{1}-\lambda _{0}}}  \label{assump0}
\end{equation}
uniformly in $k.$ From (\ref{atilde}) we have $\log \widetilde{a}\left(
n,k\right) =t_{k}-t_{k+1}$, so by (\ref{btilde}),
\begin{equation*}
\log \widetilde{b}\left( n,k\right) =\log \left( b\left( n,k\right)
e^{\lambda _{0}\left( t_{k}-t_{k+1}\right) }\right) =\log b\left( n,k\right)
+\lambda _{0}\left( t_{k}-t_{k+1}\right).
\end{equation*}

Now by assumption (\ref{assumptonb}),
\begin{equation*}
\frac{\log \widetilde{b}\left( n,k\right) }{\log \widetilde{a}\left(
n,k\right) }=\frac{\log b\left( n,k\right) +\lambda _{0}\left(
t_{k}-t_{k+1}\right) }{t_{k}-t_{k+1}}=\frac{\log b\left( n,k\right) }{%
t_{k}-t_{k+1}}+\lambda _{0}\rightarrow \lambda _{2},
\end{equation*}
finishing the proof.
\end{proof}

\section{Equidistant eigenvalues}

In this section we want to illustrate our results when the eigenvalues in $%
\left( \lambda _{0},....,\lambda _{n}\right) $ are equidistant, i.e., when $%
\lambda _{j}=\lambda _{0}+j\omega $ for $j=0,...,n.$ In this case the
elements of $E_{\left( \lambda _{0},...,\lambda _{n}\right) }$ are also
called $\mathcal{D}$-polynomials, see \cite{MoNe00} or \cite[Remark 2.1]
{GoNe}. An important particular instance of $\mathcal{D}$-polynomials is the
class of scaled trigonometric polynomials, defined for even $n$ by
\begin{equation*}
\text{span }\{1,\sin \left( 2x/n\right) ,\cos \left( 2x/n\right) ,\sin
\left( 4x/n\right) ,\cos \left( 4x/n\right), \dots,\sin x,\cos x\},
\end{equation*}
and for odd $n$ by
\begin{equation*}
\text{span }\{\sin \left( x/n\right) ,\cos \left( x/n\right) ,\sin \left(
3x/n\right) ,\cos \left( 3x/n\right), \dots,\sin x,\cos x\},
\end{equation*}
see \cite{MoNe00}.

We shall assume that $\omega \neq 0$ since $\omega =0$ leads to the
polynomial case, covered by the classical Bernstein theorem. The following
result was proved in \cite{MoNe00}:

\begin{proposition}
Let $\lambda _{j}=\lambda _{0}+j\omega $ for $j=0,...,n$, where $\omega \neq
0.$ Define $t_{k}:=a+\frac{k}{n}\left( b-a\right) $ and $p_{\left( \lambda
_{0},...,\lambda _{n}\right) ,k}$ as in (\ref{eqdefppp}). Then the operator
defined for $f\in C\left[ a,b\right] $ by
\begin{equation*}
B_{\left( \lambda _{0},...,\lambda _{n}\right) }f\left( x\right)
=\sum_{k=0}^{n}f\left( t_{k}\right) \frac{n!}{\left( n-k\right) !}\frac{%
\omega ^{k}e^{-\lambda _{0}\left( \frac{k}{n}\left( b-a\right) \right) }}{%
\left( e^{\omega \left( b-a\right) }-1\right) ^{k}}p_{\left( \lambda
_{0},...,\lambda _{n}\right) ,k}\left( x\right) 
\end{equation*}
satisfies $B_{\left( \lambda _{0},...,\lambda _{n}\right) }\left( e^{\lambda
_{0}x}\right) =e^{\lambda _{0}x}$ and $B_{\left( \lambda _{0},...,\lambda
_{n}\right) }\left( e^{\lambda _{n}x}\right) =e^{\lambda _{n}x}.$
\end{proposition}

\begin{proof}
Straightforward calculations show that the constants $d_{k}$ and $D_{k}$ in
Theorem \ref{ThmBern} are
\begin{align*}
d_{k}& :=\lim_{x\rightarrow b}\frac{\frac{d}{dx}p_{\left( \lambda
_{0},...,\lambda _{n}\right) ,k}\left( x\right) }{p_{\left( \lambda
_{1},...,\lambda _{n}\right) ,k}\left( x\right) }=-\frac{\left( n-k\right)
\omega }{1-e^{\omega \left( a-b\right) }}e^{\left( b-a\right) \left( \lambda
_{0}-\lambda _{1}\right) } \\
D_{k}& :=\lim_{x\rightarrow b}\frac{\frac{d}{dx}p_{\left( \lambda
_{0},...\lambda _{n}\right) ,k}\left( x\right) }{p_{\left( \lambda
_{0},....,\lambda _{n-1}\right) ,k}\left( x\right) }=-\frac{\left(
n-k\right) \omega }{1-e^{\omega \left( a-b\right) }}.
\end{align*}
By (\ref{eqdefft}), $t_{k}-t_{k-1}$ is defined through
\begin{equation}
e^{\left( \lambda _{0}-\lambda _{n}\right) \left( t_{k}-t_{k-1}\right) }=%
\frac{d_{k-1}}{D_{k-1}}=e^{\left( b-a\right) \left( \lambda _{0}-\lambda
_{1}\right) },  \label{defdef}
\end{equation}
so $t_{k}-t_{k-1}=\frac{\lambda _{0}-\lambda _{1}}{\lambda _{0}-\lambda _{n}}%
=\frac{1}{n}\left( b-a\right) .$ It follows that $t_{k}=a+\frac{k}{n}\left(
b-a\right) .$ According to (\ref{eqdefak}) we have
\begin{equation*}
\alpha _{k}=e^{-\lambda _{0}\left( t_{k}-a\right) }\left( -1\right)
^{k}d_{0}....d_{k-1}=e^{-\lambda _{0}\left( \frac{k}{n}\left( b-a\right)
\right) }\frac{e^{\left( a-b\right) k\omega }\omega ^{k}}{\left( 1-e^{\omega
\left( a-b\right) }\right) ^{k}}\frac{n!}{\left( n-k\right) !}.
\end{equation*}
\end{proof}

The following theorem was proved by S. Morigi and M. Neamtu in \cite[p. 137]
{MoNe00}.

\begin{theorem}
\label{ThmMN}Let $\mu _{0}\neq \mu _{1}$ be either real numbers or complex
conjugates, and in the latter case assume that $b-a<\pi /\left| Im\mu
_{0}\right| .$ Set $\Delta :=\mu _{1}-\mu _{0}$, and define $\lambda
_{j}=\mu _{0}+j\frac{1}{2n}\Delta $ for $j=0,...,2n$. Then $B_{\left(
\lambda _{0},...,\lambda _{2n}\right) }f\left( x\right) $ converges
uniformly to $f$ for all $f\in C\left[ a,b\right] .$
\end{theorem}

It is possible to derive Theorem \ref{ThmMN} from Theorem \ref{ThmCon} (for
vectors $\Lambda _{2n}$ with even index, which guarantees that $\lambda _{n}$
is a component of $\Lambda _{2n}$ for every $n)$ applied to the triple $%
\lambda _{0},\lambda _{2n}$ and $\lambda _{n}$ (instead of $\lambda
_{0},\lambda _{1},\lambda _{2})$. Since the proof is rather technical it is
omitted.

\section{Final remarks}

For $\lambda _{0}=0<\lambda _{1}<\lambda _{2}<...< \lambda _{n},$ the
M\"{u}ntz polynomials are defined as elements in the linear space $V_{n}$
generated by $1,x^{\lambda _{1}},x^{\lambda _{2}},...,x^{\lambda _{n}}$,
considered as functions on the interval $\left[ a,b\right]$, where $a\geq 0.$
Assume $a>0$. Using the transformation $x=e^{t}$ we see that $V_{n}$ is
isomorphic to the linear space $U_{n}$ generated by $1,e^{\lambda
_{1}t},e^{\lambda _{2}t},....e^{\lambda _{n}t}$. Clearly $f\in U_{n}$ is an
exponential polynomial on the interval $\left[ \ln a,\ln b\right] $. To each
$\Lambda _{n}:=\left( \lambda _{0},...,\lambda _{n}\right) $ we can
associate the Bernstein operator $B_{\left( \lambda _{0},...,\lambda
_{n}\right) }$ for the interval $\left[ \ln a,\ln b\right] $. Convergence of
$B_{n}$ to the identity operator implies that the union of the spaces $%
U_{n},n\in \mathbb{N},$ is dense in $C\left[ \ln a,\ln b\right] .$ It is
well known (see \cite[p. 180]{BoEr95}) that this entails
\begin{equation}
\sum_{n=1}^{\infty }\frac{1}{\lambda _{n}}=\infty .  \label{eqmuntz}
\end{equation}
In particular, it follows that Theorem \ref{ThmCon} does not extend to the
case of arbitrary vectors $\Lambda _{n}=\left( \lambda _{0},...,\lambda
_{n}\right)$. It would be interesting to derive from Theorem \ref{ThmCon} a
Bernstein type result for M\"{u}ntz polynomials over $\left[ \ln a,\ln b%
\right]$, $a>0$. Let us mention that the results in \cite{HiWi49} for
M\"{u}ntz polynomials over $\left[ 0,1\right] $ are of a different type,
since the basic functions used there do not form a Bernstein basis in our
sense.

\end{document}